\RequirePackage{fix-cm}

\makeatletter
\def\cl@chapter{}
\makeatother

\documentclass{svjour3}                     % onecolumn (standard format)
\smartqed  % flush right qed marks, e.g. at end of proof
\usepackage{graphicx}
\usepackage{amssymb,amsmath,amsfonts,bm}
\usepackage[margin=1in]{geometry}
\usepackage{cuted}
\usepackage[ansinew]{inputenc}
\usepackage[english]{babel}
\usepackage{fancyhdr}
\usepackage{anyfontsize}
\usepackage[active]{srcltx}
\usepackage{soul}
\usepackage{xcolor}
\usepackage[colorlinks=true, allcolors=blue]{hyperref}
\usepackage{color}
\usepackage{array}
\usepackage{textcomp}
\usepackage[linesnumbered,ruled,vlined]{algorithm2e}
\usepackage{pgfplots}
\pgfplotsset{compat=1.17}

\newcommand{\norm}[1]{\left\lVert#1\right\rVert}

\newcommand{\bfm}[1]{\mathbf{#1}}

\newcommand{\bff}{\bfm{f}}

\newcommand{\bfy}{\bfm{y}}

\newcommand{\bfz}{\bfm{z}}

\newcommand{\alphavec}{\boldsymbol \alpha}

\newcommand{\bfM}{\bfm{M}}

\newcommand{\bfH}{\bfm{H}}

\newcommand{\bfA}{\bfm{A}}

\newcommand{\bfB}{\bfm{B}}

\newcommand{\bfE}{\bfm{E}}
\newcommand{\bfF}{\bfm{F}}

\newcommand{\bfI}{\bfm{I}}
\newcommand{\bfX}{\bfm{X}}

\newcommand{\bfY}{\bfm{Y}}

\newcommand{\bfQ}{\bfm{Q}}

\newcommand{\bfR}{\bfm{R}}

\newcommand{\bfW}{\bfm{W}}

\newcommand{\bfphi}{\pmb{\phi}}

\usepackage{subfiles}
\usepackage{subcaption}
\usepackage{smartdiagram}
\usesmartdiagramlibrary{additions}
\usepackage{lineno}
\usepackage{upgreek}
\usepackage{textgreek}
\usepackage{cancel}
\usepackage[capitalise]{cleveref}
\Crefname{equation}{Eq.}{Eqs.}
\Crefname{appendix}{Appendix}{Appendices}
\usepackage{comment}
\usepackage{svg}
\usepackage{float}
\allowdisplaybreaks
\usepackage{nicefrac}
\graphicspath{ {./figures/} }
\usepackage[export]{adjustbox}
\usepackage{numprint}

\modulolinenumbers[5]

\journalname{}
\begin{document}

\title{On the validity limits of the parametrisation method for invariant manifolds}
\subtitle{An assessment of practical criteria for vibrating systems}
\titlerunning{On the validity limits of the parametrisation method for invariant manifolds}        % if too long for running head

\author{A. de Figueiredo Stabile \and A. Grolet \and A. Vizzaccaro \and C. Touz{\'e}}


\institute{A. de Figueiredo Stabile \at
           Institute of Mechanical Sciences and Industrial Applications (IMSIA), ENSTA - CNRS - EDF, Institut Polytechnique de Paris, Palaiseau, France. \\
           \email{andre.de-figueiredo-stabile@ensta.fr}
           \and
           A. Grolet \at
           LISPEN, Arts et Metiers Sciences et Technologies, Lille, France.
           \and
           A. Vizzaccaro \at
           College of Engineering, Mathematics and Physical Sciences, University of Exeter, Exeter, United Kingdom.
           \and
           C. Touz{\'e} \at
           Institute of Mechanical Sciences and Industrial Applications (IMSIA), ENSTA - CNRS - EDF, Institut Polytechnique de Paris, Palaiseau, France.
}

\date{}

\maketitle

\begin{abstract}
The parametrisation method for invariant manifolds is a powerful technique for deriving reduced-order models in the context of nonlinear vibrating systems, allowing accurate computations of nonlinear normal modes. Thanks to arbitrary order asymptotic expansions, converged results are within reach and directly applicable to finite element structures. However, since it relies on a local theory and asymptotic expansions, the results are only valid up to a given amplitude, which defines the convergence radius of the approximation. The aim of this contribution is to investigate the validity limits of the approach and review the existing error estimates, with the concrete objective of proposing a practical approach to estimate the validity range during the computation, thus producing safe bounds within which the reduced-order model can be used. Three different criteria are assessed. The first one uses the error in the invariance equation as the distance to the fixed point increases. The second one is adapted from an upper bound criterion derived for normal form transforms and based on the potential singularities of the homological operator. The third one uses validity limits of series expansion through Cauchy and d'Alembert rules, which can be tested either on the reduced dynamics coefficients or those of the nonlinear mappings. The criteria are tested on a number of different examples that are representative of the situations encountered when dealing with nonlinear vibrations. The Duffing equation serves as a first benchmark that allows considering conservative oscillations, forced systems at primary resonance, and superharmonic resonance. The investigations are then extended to a vibrating system with two degrees of freedom. Finally, the different criteria are assessed on a finite element beam structure, and guidelines are formulated to generalise their practical use and produce accurate and easy-to-use error bounds in the context of model order reduction for nonlinear vibrating structures.
\end{abstract}

\section{Introduction}

The parametrisation method for invariant manifolds is an efficient technique which is used in dynamical systems theory to prove existence and uniqueness of invariant manifolds~\cite{Cabre1,Cabre2,Cabre3,Haro,REINHARDT2019,VANDENBERG:center}, but also as an efficient algorithm that produces high-order approximations which can be employed for simulation-free model order reduction. In particular, it has been used successfully in the field of nonlinear vibrations in recent years, proposing accurate and effective reduced-order models (ROMs) by leveraging the idea of Nonlinear Normal Modes (NNMs) as invariant manifolds of the phase space~\cite{PONSIOEN2018,PONSIOEN2020,JAIN2021How,vizza21high,opreni22high,Vizza:superDPIM}. Even though elements of the method were already present at various points in earlier developments, see {\em e.g.}~\cite{Coullet83,Roberts89,RobertsRod93,Carini2015}, only the systematic procedure as summarized {\em e.g.} in~\cite{Haro} provides a clear understanding of the existence of an infinite number of parametrisation styles, among which two are of particular significance: the graph style and the normal form style. In the field of nonlinear vibrations, it allows unifying different techniques proposed to compute NNMs for model-order reduction: the center manifold technique \cite{ShawPierre91,ShawPierre93,ShawShock} and the normal form approach~\cite{Jezequel91,touze03-NNM,LEUNG2003,TOUZE:JSV:2006,ReviewROMGEOMNL,TouzeCISM2}. In recent years, an important step forward has been taken by proposing direct versions of the parametrisation method, in the sense that the technique can be directly applied to the physical space, and only the computation of the master eigenmodes is necessary. In the remainder of the article, we will refer to the DPIM as Direct Parametrisation for Invariant Manifolds, for all the methods proposing such a treatment~\cite{JAIN2021How,vizza21high}. Thanks to that,  applications to industrial-scale models have been reported and importantly extend the applicability range of the method~\cite{vizza21high,opreni22high,li2021periodic,Martin:rotation,MingwuLi2024,frangi_electromechDPIM,Pinho:shells,Stabile:follow}.

Although the parametrisation method is highly versatile and allows for accurate results thanks to high-order asymptotic expansions, it suffers from the same inconvenience as any local approach~\cite{Haro,gucken83}: once a certain convergence radius, or validity limit, of the expansion is reached, the results diverge. A few strategies have been proposed to deal with this limitation, such as using Pad{\'e} approximants to patch a global approximation to the local one~\cite{Kaszas2025}, or by combining a restart procedure obtained from the asymptotic-numerical method (ANM) with a harmonic balance normal form parametrisation~\cite{Grolet:HBNF}. Nevertheless, determining the convergence domain of the approximation remains an important point to address. In~\cite{Haro}, an \emph{a posteriori} error estimate is proposed in terms of the fulfillment of the invariance equation. Such a criterion allows defining the amplitude limit over which the approximations fail to stick to the true invariant manifold, and has been used in different mathematical contexts~\cite{Castelli2015,Breden2016,VanDenBerg2016,Groothedde:dde}. In the field of nonlinear vibrations, the use of classical series convergence rules (d'Alembert and Cauchy) has also been proposed to quantify the maximum amplitude up to which convergence of the series can be ensured~\cite{Grolet2025}. When restricting to the normal form parametrisation, the upper bound derived in~\cite{LamarqueUP}, which is based on the singularity of the homological operator, might also be used, even though asking for an adjustment to handle the more generic case of the parametrisation method. Although the validity limits of the method are a critical issue, recent literature on its applications to vibrating structures seldom addresses them. Even if in the convergence range the asymptotic expansions converge and can thus be safely employed to propose quasi-exact solutions, the analyst needs to have an accurate idea of the maximum amplitude limit up to which the method can be safely used.

The aim of this contribution is thus to review the existing criteria for determining the validity limit of the parametrisation method, and compare them in the specific framework of nonlinear vibrating systems. Also, an adaptation of the upper bound introduced in~\cite{LamarqueUP} is proposed to fit with the DPIM. More specifically, the goal is to assess the behaviour of the different criteria for the cases at hand in nonlinear vibrations. To that purpose, free oscillations of conservative systems will be studied, but also forced oscillations (non-autonomous systems) in primary and secondary resonances scenarios. Only the reduction to a single NNM is considered here, while the case of multiple modes is left for future work. After this assessment, a practical methodology is proposed for the analyst, in order to give the user accurate yet easy-to-use bounds, which can predict beforehand the amplitude range in which ROMs based on the DPIM can be safely used.

The text is divided as follows: \cref{sec:ValLims} presents the main theoretical developments, with a brief summary of the considered algorithm for the determination of the DPIM, followed by a presentation of the three different criteria and some details on their practical implementation. Next, in \cref{sec:Examples}, several different applications are presented. First, 
the case of the Duffing equation under various circumstances is investigated. In this scenario, no reduction is performed, and some simplifications occur when applying the criteria. Then, a discrete system with two degrees of freedom (DOF) is considered, where the effect of model reduction can be showcased in a small-scale problem. Finally, a beam discretised by finite elements is studied, to exemplify the results on a large-scale model. Once all of the cases have been investigated, a methodology is proposed to obtain effective and rapid computations of the range of convergence from a practical viewpoint, addressing in particular issues pertaining to their use when only low-order computations are available.

\section{Validity limits of asymptotic expansions in the context of the parametrisation method} \label{sec:ValLims}

The aim of this section is to introduce the different criteria that will be assessed for predicting the validity limit of asymptotic expansions derived in the context of the parametrisation method for invariant manifolds applied to nonlinear vibrating systems. The method is first briefly recalled for self-consistency.

\subsection{The parametrisation method for invariant manifolds}

\subsubsection{Autonomous systems}

The direct parametrisation method (DPIM) for autonomous systems is first reviewed, then the case of non-autonomous problems is recalled. Whereas descriptions of the parametrisation method in its general settings can be found {\em e.g.} in~\cite{Haro}, the focus is here set to its application to nonlinear vibrating systems. The method is also termed as {\em direct} as it can be applied directly to the physical space. Different technical adjustments have been proposed to make this possible, using {\em e.g.} a norm-minimizing method~\cite{JAIN2021How}, or a bordering technique~\cite{Carini2015,vizza21high}. The presentation retained here follows the guidelines given in~\cite{vizza21high}, and subsequently expanded in \cite{Vizza:superDPIM} to include forcing terms of arbitrary order.

To make the presentation lighter, the general form of the problem at hand is considered to be written as an autonomous quadratic system of differential algebraic equations (DAEs):
\begin{equation} \label{eq:autonomous_EOM}
	\bfB \dot{\bfy} = \bfA \bfy + \bfQ(\bfy,\bfy).
\end{equation}
In the above, $\bfy \in \mathbb{C}^D$ is the state vector, $\bfA,\bfB \in \mathbb{R}^{D\times D}$ are matrices, with $\bfA$ non-singular, and $\bfQ:\mathbb{R}^D \times \mathbb{R}^D \to \mathbb{R}^D$ is a quadratic operator responsible for the nonlinear forces.

\begin{remark}
An important aspect granting flexibility to the DPIM is the fact that the matrix $\bfB$ can be singular. This property endows the method with the ability to treat systems of DAEs automatically, allowing in particular the use of quadratic recast to deal with arbitrary polynomial nonlinearities \cite{Guillot:recast,Grolet2025}. This means that the option to consider only the quadratic operator as a source of nonlinear effects is not restrictive, as any type of analytical nonlinearity can be recast into a quadratic form~\cite{KARKAR2013,Guillot:recast,QuadRecast}. It is merely a choice made in order to allow for an easier and unified treatment of the nonlinearities. When quadratic recast is employed, supplementary \textit{algebraic} equations, with left-hand side equal to zero, defining auxiliary variables needed for its use, are appended to the original system. If this is the case, in all of the developments presented in this contribution, they will be placed in the last $D_a$ entries of the state vector, and the first block of dimension $D_p \times D_p$, with $D = D_p + D_a$, of the matrix $\bfB$, related to the physical variables, will be supposed non-singular.
\end{remark}

The method's main idea is, then, to look for invariant manifolds of a fixed point of the system in \cref{eq:autonomous_EOM}. To this end, a coordinate system of so-called normal coordinates $\bfz$ is introduced, and a parametrisation $\bfW(\bfz)$ of the manifold representing the physical coordinates is sought in the form of a polynomial expansion in $\bfz$. The method additionally gives an evolution equation, called the reduced dynamics and denoted $\bff(\bfz)$, for the normal coordinates on the manifold. The series expansions for these quantities are then written as
\begin{subequations} \label{eq:series_exp}
	\begin{align} 
		\bfy &=  \bfW(\bfz) = \sum_{p=1}^{o} \sum_{k=1}^{m_p} \bfW^{(p,k)} \bfz^{\alphavec(p,k)}, \\
		\dot{\bfz} &=\bff(\bfz) = \sum_{p=1}^{o} \sum_{k=1}^{m_p} \bff^{(p,k)} \bfz^{\alphavec(p,k)},
	\end{align}
\end{subequations}
with $o$ the maximum order of the parametrisation, $m_p= \begin{pmatrix} p+d-1 \\ p \end{pmatrix}$ the number of order $p$ monomials, $\bff^{(p,k)},\bfW^{(p,k)}$ the coefficients of the polynomial expansions, $\bfz$ the normal coordinates vector of dimension $d$, and $\alphavec(p,k)$ the $k-$th multi-index vector of order $p$, such that $\alphavec(p,k)_i \in \mathbb{N}$, $|\alphavec(p,k)|=p$, and $\bfz^{\alphavec(p,k)} = z_1^{\alphavec(p,k)_1} \cdot \ldots \cdot z_d^{\alphavec(p,k)_d}$. Note that $\alphavec(p,k)$ is a vector of integers and $\alphavec(p,k)_i$ denotes the $i$-th entry.

\begin{remark}
It is here assumed that the fixed point of the system is located at the origin of the state space. The case of non-zero fixed points only entails the use of tangent operators, but does not substantially change the analysis procedure. In such a case, a simple shift can be operated to relocate the fixed point at the origin, see {\em e.g.}~\cite{Martin:rotation,frangi_electromechDPIM,Stabile:follow,Grolet2025}.
\end{remark}

\begin{remark}
The dimension of the normal variable vector is typically much smaller than that of the state vector, {\em i.e.} $d \ll D$, which makes the method an efficient reduction technique. The choice of the dimension $d$ of the manifold is related to {\em a priori} knowledge of the system~\cite{ReviewROMGEOMNL}. In particular, for lightly damped mechanical systems, $d$ is always even, and its precise value is determined by the potential occurrence of \textit{internal resonances}, {\em i.e.} commensurability relations between the linear undamped eigenfrequencies of the system.
\end{remark}

By plugging \cref{eq:series_exp} into \cref{eq:autonomous_EOM} and using the chain rule, the time dependence of \cref{eq:autonomous_EOM} can be discarded, giving rise to the invariance equation:
\begin{equation} \label{eq:invariance}
	\bfB \nabla_{\bfz} \bfW(\bfz) \bff(\bfz) = \bfA \bfW(\bfz) + \bfQ(\bfW(\bfz),\bfW(\bfz)).
\end{equation}
By isolating terms of different polynomial orders, the invariance equation can be split into several homological equations, amenable to an order-by-order solution. In particular, at each order, these equations decouple for the different monomials on the expansions given by \cref{eq:series_exp}, resulting in a series of linear systems to be solved. Using a bordering technique which augments the size of the ill-conditioned systems to be solved by adding an orthogonality condition to the left eigenvectors corresponding to its kernel, the system for the $k-$th monomial of order $p$ finally reads~\cite{Vizza:superDPIM}:
\begin{equation} \label{eq:homological}
	\begin{bmatrix}
		\sigma^{(p,k)} \bfB - \bfA & \bfB \bfY_{\mathcal{R}} & \bf0
		\\
		\bfX_{\mathcal{R}}^\star \bfB &  \bf0 & \bf0 \\
		\bf0 & \bf0 & \bfI
	\end{bmatrix}
	\begin{bmatrix}
		\bfW^{(p,k)}\\
		\bff^{(p,k)}_{\mathcal{R}}\\
		\bff^{(p,k)}_{\cancel{\mathcal{R}}}
	\end{bmatrix}
	=
	\begin{bmatrix}
		\bfR^{(p,k)}\\
		\bf0 \\
		\bf0
	\end{bmatrix}.
\end{equation}
In the above equation, $\bfX$ and $\bfY$ are matrices grouping the left and right master eigenvectors of the system, and $\bfR^{(p,k)}$ denotes the right-hand side, or residual, vector for the considered monomial, and depends only on previously computed quantities. Its explicit expression can be found in \cite{Vizza:superDPIM}, with an extension in \cite{Stabile:follow} allowing to treat non-diagonal reduced dynamics. Additionally, the definition
\begin{equation}
	\sigma^{(p,k)} = \sum_{i = 1}^d \alphavec(p,k)_i \cdot \lambda_i
\end{equation}
holds, with $\lambda_i$ the eigenvalues of the system. The subscript $\mathcal{R}$ allows to select only certain eigenvalues from the variables in which it appears, and is related to the so-called resonance conditions of the system: if $s \in \left[ 1, \ldots, d \right]$ is in $\mathcal{R}$ (it is then not in $\cancel{\mathcal{R}}$), the monomial $(p,k)$ is not set to zero in the reduced dynamics equation for variable $z_s$. When $\sigma^{(p,k)} \approx \lambda_s$ for some $s \in \left[ 1, \ldots, d \right]$, a resonance is said to occur, and the index $s$ is necessarily kept into~$\mathcal{R}$. At this point, a few choices can be made regarding different \textit{styles} of parametrisation. If only the monomials for which $\sigma^{(p,k)} \approx \lambda_s$ are kept in the reduced dynamics, then the complex normal form (CNF) style is selected. Instead, if the definition of $\mathcal{R}$ is enlarged to its maximum, such that no monomial is set to zero in the reduced dynamics - which corresponds to the set of non-resonant monomials $\cancel{\mathcal{R}}$ being empty and the last line of \cref{eq:homological} vanishing - the resulting parametrisation is called the graph style. Other than these two canonical forms, an infinity of other styles is also possible, relating to different choices on how to populate the set $\mathcal{R}$: real normal form, oscillator normal form and mixed styles, see {\em e.g.} \cite{Haro,vizza21high,opreni22high,TouzeCISM2}.

\begin{remark}
When \cref{eq:homological} is considered for $p = 1$, it yields the right eigenvectors as the linear coefficients of the parametrisation, {\em i.e.} $\bfW^{(1,k)} = \bfY_k$, and a diagonal linear reduced dynamics provided that no Jordan blocks are present. This means that the searched invariant manifolds are tangent to the linear eigenspaces corresponding to the master eigenvectors $\bfY_k$. In particular, the choice of how many (and which) eigenvectors to include in the linear part will dictate which invariant manifold is parametrised. Usually, for a lightly damped mechanical system, a single NNM, associated with a complex conjugate pair of eigenvalues, can be selected, resulting in a two-dimensional invariant manifold. This choice of a single NNM reduction will be retained in the present study, but other choices are possible, {\em e.g.} for parameter-dependent systems experiencing bifurcations \cite{Stabile:follow} or the study of internal resonances \cite{AndreaROM,opreni22high,Vizza:superDPIM,li2021periodic,li2021bifurcation}.
\end{remark}
	
\begin{remark}
Although the parametrisation method is most frequently used for model-order reduction, when $d \ll D$, it can also be used for the case where no reduction takes place, {\em i.e.} $d=D$. In this scenario, if a CNF style is chosen, the method corresponds to the usual normal form technique~\cite{Haro,Jezequel91,AndreaROM}, and allows for a simplification of the problem dynamics, keeping only the resonant monomials.
\end{remark}

\subsubsection{Non-autonomous systems}

The case of non-autonomous systems can be seen as an extension of the autonomous case by suitably choosing additional normal coordinates in order to account for the forcing~\cite{Vizza:superDPIM}. We consider, in particular, the case of a single harmonic excitation, but an extension for multiple harmonics is readily available. We consider a system of the form
\begin{equation} \label{eq:nonaut_syst}
	\bar{\bfB} \dot{\bar{\bfy}} = \bar{\bfA} \bar{\bfy} + \bar{\bfQ}(\bar{\bfy},\bar{\bfy}) + \bfF_c \cos{\Omega t} + \bfF_s \sin{\Omega t},
\end{equation}
where $\bfF_c$ and $\bfF_s$ are spatial distributions of the forces and the introduction of the barred notation will become clear shortly. The forcing can be written in complex form as
\begin{equation}
	\bfF_c \cos{\Omega t} + \bfF_s \sin{\Omega t} = \tilde{\bfF} \tilde{\bfy},
\end{equation}
with 
\begin{equation}
	\tilde{y}_{1,2} = e^{\pm i \Omega t},
\end{equation}
and the columns of $\tilde{\bfF}$ suitably chosen. Then, by noticing that $\dot{\tilde{y}}_{1,2} = \pm i \Omega y_{1,2}$, the vector $\tilde{\bfy}$ can be seen as composed of new auxiliary coordinates, and the forced system can be recast as an autonomous one:
\begin{equation} \label{eq:nonautonomous_EOM}
	\underbrace{
		\begin{bmatrix}
			\bar{\bfB} & \bf0 \\
			\bf0 & \tilde{\bfB}
		\end{bmatrix}
	}_{\bfB}
	\dot{
		\underbrace{
			\begin{bmatrix}
				\bar{\bfy} \\
				\tilde{\bfy}
			\end{bmatrix}
		}_{\dot{\bfy}}
	}
	= 
	\underbrace{
		\begin{bmatrix}
			\bar{\bfA} & \tilde{\bfF} \\
			\bf0 & \tilde{\bfA}
		\end{bmatrix}
	}_{\bfA}
	\underbrace{
		\begin{bmatrix}
			\bar{\bfy} \\
			\tilde{\bfy}
		\end{bmatrix}
	}_{\bfy}
	+
	\underbrace{
		\begin{bmatrix}
			\bar{\bfQ}(\bar{\bfy},\bar{\bfy}) \\
			\bf0
		\end{bmatrix}
	}_{\bfQ(\bfy,\bfy)},
\end{equation}
with $\bfy$ now of dimension $D+2$, and matrices $\tilde{\bfA}$ and $\tilde{\bfB}$ defined as
\begin{equation} \label{eq:AandBforcing}
	\tilde{\bfB} = 
	\begin{bmatrix}
		1 & 0 \\
		0 & 1
	\end{bmatrix},
	\quad
	\tilde{\bfA} = 
	\begin{bmatrix}
		i \Omega & 0 \\
		0 & -i \Omega
	\end{bmatrix}.
\end{equation}
The parametrisation method can then be applied to \cref{eq:nonautonomous_EOM} in the same fashion as described before, with only minor modifications. Now, the vector $\bfz = \left[ \bar{\bfz}^T \quad \tilde{\bfz}^T \right]^T$ has dimension $d+2$, and the coordinates $\tilde{\bfy}$ are assumed to already be in normal form, {\em i.e.} $\tilde{\bfy} = \tilde{\bfz}$. Therefore, the dynamics for $\tilde{\bfz}$ is already known, and both the parametrisation and the reduced dynamics will only be searched for the autonomous coordinates, but including monomials also with the non-autonomous ones. For the interested reader, further details on this procedure, as well as a comparison to approaches where the forcing is supposed to be a small perturbation, can be found in~\cite{Vizza:superDPIM}.

We now turn to the presentations of the three different available criteria which will be assessed. A specific emphasis will be put on the case of the reduction to a single NNM (two-dimensional invariant manifold), since this case will be studied in Section~\ref{sec:Examples}.

\subsection{Error in the invariance equation}

The first considered validity criterion is related to the invariance equation,~\cref{eq:invariance}. It has been developed in the mathematical community, and is widely used in all their examples in order to assess the quality of the ROM given by the method, see~{\it e.g.}~\cite{Haro,Breden2016,VanDenBerg2016,REINHARDT2019}. It relies on a quantitative measure of the error in the invariance of the manifold itself, obtained by checking how well the invariance property is satisfied by the asymptotic expansion, which may diverge from the exact solution.
As the solution to \cref{eq:invariance} is found only in an approximate manner, it is possible to consider that the approximation fails once the norm of the difference between the left-hand and right-hand sides of the equation is larger than a certain threshold. Specifically, we take
\begin{equation} \label{eq:inv_eq_error}
	\frac{\norm{\bfB \nabla_{\bfz} \bfW(\bfz) \bff(\bfz) - \bfA \bfW(\bfz) - \bfQ(\bfW(\bfz),\bfW(\bfz))}}{\norm{\bfA \bfW^{(1,1)}}} < \varepsilon
\end{equation}
to define the validity limit of the approximation. In the above equation, the factor $\norm{\bfA \bfW^{(1,1)}}$ is introduced to normalise the residual of the invariance equation, which is dimensionally homogeneous to a force in the case of vibrating systems. A drawback of this approach is the need to arbitrarily select a tolerance value.

In practice, assuming that reduction to a single NNM is undertaken and that the problem is autonomous, the computed invariant manifold is two-dimensional and described with the normal variables $\bfz=[z_1 \; z_2]^T$, which are complex conjugate, {\em i.e.} $z_2=\bar{z}_1$. The computation of the criterion given by~\cref{eq:inv_eq_error} thus needs to be performed in a two-dimensional domain. Since our main objective is to compute a radius of convergence, the use of a polar representation for the complex coordinates is natural. Thus, considering $\rho$ as a measure of their amplitude and the phase $\theta$, one writes
\begin{equation} \label{eq:polar_autonomous}
	z_{1,2} = \frac{\rho}{2} e^{\pm i \theta}.
\end{equation}
The calculation then proceeds by discretising the interval $[0,2\pi]$ for $\theta$ and applying a root-finding algorithm to determine a validity limit for each of the resulting angles. The global validity limit is then computed in the end as the minimum $\rho$ obtained for the several different sampled angles. This approach is particularly appealing, as it enables the use of single-variable root-finding algorithms, which are far less computationally expensive and considerably easier to implement than their multidimensional counterparts.

If we now consider the case of a non-autonomous system, still parametrised by a single master mode, the auxiliary coordinates related to the forcing can be written as
\begin{equation} \label{eq:polar_non_autonomous}
	z_{3,4} = e^{\pm i \Omega t},
\end{equation}
with $\Omega$ the forcing frequency. In this situation, by further assuming the somewhat general case of an $m:n$ resonance, the response is periodic, and such that
\begin{equation} \label{eq:resonance}
	\Omega = \frac{m}{n} \dot{\theta} \Rightarrow \Omega t = \frac{m}{n} \theta - \phi,
\end{equation}
where $\phi$ is a phase which represents the shift between the forcing and the response. The introduction of the forcing, then, implies considering an additional coordinate in the domain when looking for the convergence radius. The same procedure as for the autonomous case is followed: several different values of $\theta$ and $\phi$ are computed, and the minimum value of $\rho$ obtained for all of them is taken as the validity limit of the approximation.

\begin{remark}
The angles $\theta$ and $\phi$ should be allowed to vary in intervals corresponding to the full periodic behaviour of the solution. This might correspond to the interval $\left[0, 2\pi \right]$ for both angles or to larger intervals for the angle $\theta$, in the case of secondary resonances ({\em e.g.} for a 1:3 superharmonic resonance the chosen interval must be $\left[0, 6\pi \right]$).
\end{remark}

Interestingly, a more direct version of the criterion can also be employed, by using an idea devised in the context of the asymptotic numerical method (ANM), which consists in assuming that most of the error in asymptotic expansions is contained within the first neglected term: if the asymptotic expansion goes up to a given order $o$, then the term of order $o+1$ contains most of the error as leading-order residual~\cite{Cochelin1994}. With such an idea and using polar representation in the single NNM autonomous case, an explicit expression for the maximum amplitude can then be derived, making the calculation of the validity limit straightforward. In this section, only the main result and a brief summary of its derivation are given. The reader is referred to~\cref{sec:inv_eq_simpl} for a detailed description. 

The calculation of the error can be used in the present context by combining the residual vectors $\bfR^{(p,k)}$ appearing in~\cref{eq:homological}, into a single vector for order $p$, gathering them sequentially. The residual vector $\bfR^{(p)}$ at order $p$ is then defined as
\begin{equation}
	\bfR^{(p)} = \left[ (\bfR^{(p,k_1)})^T \quad (\bfR^{(p,k_2)})^T \quad \ldots \quad (\bfR^{(p,m_p-1)})^T \quad (\bfR^{(p,m_p)})^T \right]^T.
\end{equation}
Let us denote as $\bfE$ the error in the invariance equation, given in~\cref{eq:inv_eq_error}:
\begin{equation}\label{eq:numinvaerror}
	\bfE=\bfB \nabla_{\bfz} \bfW(\bfz) \bff(\bfz) - \bfA \bfW(\bfz) - \bfQ(\bfW(\bfz),\bfW(\bfz)).
\end{equation}
As shown in~\cref{sec:inv_eq_simpl}, in the case of a computation led up to order $o$, a direct relationship exists between $\bfE$, and the residual vector $\bfR^{(o+1)}$. Consequently, a meaningful approximation consists in simplifying the full calculation of $\norm{\bfE}$ by the first neglected term in the residual, multiplied by the related power of the normal variable amplitude, as
\begin{equation} \label{eq:norm_approx_HBNF}
	\norm{\bfE}  \approx \norm{\bfR^{(o+1)}}  \left( \frac{\rho}{2} \right)^{o+1}.
\end{equation}
The above expression yields the following direct relationship between the radius of convergence $\rho$ and the tolerance $\varepsilon$:
\begin{equation} \label{eq:inv_eq_error_ANM}
	\rho = 2 \left(\frac{\norm{\bfA \bfW^{(1,1)}}}{\norm{\bfR^{o+1}}} \, \varepsilon \right)^{\frac{1}{o+1}}.
\end{equation}
An analytical expression as such cannot be obtained in the forced case because the presence of the non-autonomous coordinates creates monomials with several different powers of $\rho$ for a fixed polynomial order, thus making it impossible to isolate it straightforwardly. It would be possible, however, to determine a polynomial in $\rho$ that could be solved numerically to yield the desired convergence radius.

\subsection{Singularity of the homological operator} \label{sec:Singularity}

A criterion for determining the validity limit of normal form transforms for autonomous systems has been proposed in~\cite{LamarqueUP}. Since the parametrisation method shares common ideas with the normal form technique, especially when a CNF style is chosen, the proposed approach will be adapted here to the present context and to the case of non-autonomous systems. For the moment, we restrict our attention to situations where $d=D$, so that no reduction takes place, in order to introduce the criterion. By inspecting the invariance equation, \cref{eq:invariance}, we see that to find the reduced dynamics $\bff(\bfz)$, one needs to be able to, at least in principle, invert the homological operator $\bfH$, defined as 
\begin{equation}
	\bfH = \bfB \nabla_{\bfz} \bfW(\bfz).
\end{equation}

This inversion does not actually take place in practice, as it is too costly, and it is indeed not effectuated due to the order by order solution of the invariance equation, but it suggests that the singularities of the operator $\bfH$ bound the region at which the parametrisation remains valid. In such way, one can look for the singularities of $\bfH$ in order to determine a criterion allowing to estimate the validity limit of the approximation. For the case $d=D$ this corresponds to computing the determinant of $\bfH$, which is a square matrix, and finding its zeros. 
In~\cite{LamarqueUP}, it has been found on numerical examples that this criterion provides an effective upper bound for validity limits.

\begin{remark}
When considering DAEs, the operator $\bfH$ is automatically singular for all values of $\bfz$ because of the lines related to the algebraic equations, and this fact hinders the application of the singularity criterion directly. To circumvent this, for the case of DAEs, only the equations of the system related to the physical variables, for which the left-hand side matrix is non-singular, should be considered.
\end{remark}

Considering the general non-autonomous case, one can then split the coefficients of the nonlinear mapping $\bfW$ in order to isolate the non-autonomous coordinates $\tilde{\bfz}$ as
\begin{equation}
	\bfW(\bfz) = 
	\begin{bmatrix}
		\bar{\bfW}(\bfz) \\
		\tilde{\bfz}
	\end{bmatrix}.
\end{equation}
The homological operator can then be decomposed into several blocks:
\begin{equation}
	\bfH =
	\begin{bmatrix}
		\bar{\bfB} \nabla_{\bar{\bfz}} \bar{\bfW}(\bfz) & \bar{\bfB} \nabla_{\tilde{\bfz}} \bar{\bfW}(\bfz) \\
		\bf0 & \tilde{\bfB}
	\end{bmatrix}.
\end{equation}
Using~\cref{eq:AandBforcing}, the determinant of the homological operator then simply reads:
\begin{equation}
	\det{\bfH} = \det{\bar{\bfB} \nabla_{\bar{\bfz}} \bar{\bfW}(\bfz)},
\end{equation}
which shows that only the autonomous part of the mapping plays a role when determining the singularities of $\bfH$.

Also in this case, the same type of procedure employed for the invariance equation criterion is used: discrete values of the angle $\theta$ (and possibly $\phi$) are chosen, and the radius of convergence is determined as the minimum $\rho$ for all of them, computed each time by a root-finding algorithm.

Considering now the case where $d<D$, the homological operator is represented by a rectangular matrix, and thus a direct computation of its determinant is not possible. A conceivable alternative would be to numerically determine its rank by computing its singular value decomposition (SVD), and find the points at which it becomes smaller than $d$. This can be done by noting that non-zero singular values are always positive, and thus tracking the smallest one and employing a minimisation algorithm to determine where it vanishes. Another simpler possibility consists of selecting $d$ representative degrees of freedom for the problem at hand and considering the submatrix obtained by extracting their corresponding lines in the homological operator. The resulting matrix is square, and thus its determinant can be easily computed. This was the choice retained in this contribution, as its implementation is straightforward from the case $d=D$. 
A comparison with the singular values method showed that both approaches yield identical results. The selected test case is presented in~\cref{sec:sing_comp}. It should be noted that the singular values approach can be seen as more rigorous, and might be particularly interesting when there is no obvious choice of representative degrees of freedom to extract the relevant submatrix.
	
\subsection{Series expansions: d'Alembert and Cauchy rules}

If a single master mode is selected, the primary results of the method, {\em i.e.} nonlinear mappings and reduced dynamics, can be expressed as power series on the amplitude of the complex conjugate normal coordinates $z_{1}$ and $z_{2}$. It is then possible, as proposed in \cite{Grolet2025}, to apply traditional convergence criteria for series to assess the radius of convergence of the approximation. Since only one master mode is chosen for the parametrisation, the two (complex conjugate) normal coordinates have the same modulus, and can be written in polar form as given by \cref{eq:polar_autonomous}, with additionally the forcing being determined by \cref{eq:polar_non_autonomous}. Replacing these expressions into the polynomial expansions given by \cref{eq:series_exp}, they can be regarded as truncated power series in $\rho$:
\begin{subequations}
	\begin{align}
		\bff(\bfz) &= \sum_{p=1}^{o} \bff_p(\theta, \Omega t) \, \rho^p, \\
		\bfW(\bfz) &= \sum_{p=1}^{o} \bfW_p(\theta, \Omega t) \, \rho^p,
	\end{align}
\end{subequations}
in which the coefficients only depend on $\Omega t$ for the non-autonomous case. In this situation, by further assuming once again an $m:n$ resonance for the forcing, \cref{eq:resonance} is valid, and, by recognizing that in \cref{eq:series_exp} the exponent of variable $z_i$ for the $k-$th monomial of order $p$ is given by $\alphavec(p,k)_i$, the coefficients above can be rewritten as a function of the phase $\phi$ as
\begin{subequations}
	\begin{align}
		\bff_p(\theta, \phi) &= \sum_{k \in \mathcal{K}_p} \frac{\bff^{(p,k)}}{2^p} e^{i [\alphavec_1-\alphavec_2 + \frac{m}{n}(\alphavec_3-\alphavec_4)] \theta} e^{i(\alphavec_4-\alphavec_3)\phi} , \\
		\bfW_p(\theta, \phi) &= \sum_{k \in \mathcal{K}_p} \frac{\bfW^{(p,k)}}{2^p} e^{i [\alphavec_1-\alphavec_2 + \frac{m}{n}(\alphavec_3-\alphavec_4)] \theta} e^{i(\alphavec_4-\alphavec_3)\phi},
	\end{align}
\end{subequations}
where $\alphavec(p,k)$ has been denoted by $\alphavec$ to lighten the notation and the sets $\mathcal{K}_p$ are defined, by noting that for a given monomial the power of $\rho$ is given by $\alphavec_1+\alphavec_2$, as
\begin{equation}
	\mathcal{K}_p = \{ k \in \{1,\ldots,m_p\} : \alphavec(p,k)_1 + \alphavec(p,k)_2 = p \}.
\end{equation}

Then, for fixed $\theta$ and $\phi$ values, one can estimate a validity limit for the approximation by considering the radius of convergence of the power series in $\rho$. In order to do so, d'Alembert's and Cauchy's series convergence criteria are considered~\cite{Rudin76}, such that the estimated convergence radii can be computed by choosing a scalar component of either $\bff_p(\theta, \phi)$ or $\bfW_p(\theta, \phi)$, denoted hereafter by $a_p$, and computing the limits
\begin{align}
	\rho_C &= \lim_{p \to \infty} a_p^{-1/d}, \\
	\rho_A &= \lim_{p \to \infty} \frac{a_{p}}{a_{p+1}}.
\end{align}

Evidently these expressions cannot be exactly evaluated numerically, but provided a high enough order of parametrisation $o$ is chosen, a good approximation should be obtained. The reasoning above is valid for fixed $\theta$ and $\phi$. In order to account for the dependence on these variables the same procedure described for the other two criteria is employed, and a global radius of convergence can be determined by considering the minimum of the obtained values for all combinations of angles.

\begin{remark}
It should be noted that, often times, the series stemming from the parametrisation method only contains even or odd terms. In this case, the expression for d'Alembert's convergence radius has to be adapted to
	\begin{equation}
		\rho_A^2 = \lim_{p \to \infty} \frac{a_{2p}}{a_{2p+2}}.
\end{equation}
\end{remark}

\begin{remark}
When computational constraints limit the maximum order of parametrisation, for example in the case of high-dimensional systems discretized by finite elements, extrapolation procedures can be employed in order to determine the convergence radius from low-order computations, as will be exemplified in \cref{sec:Duffing_consevative,sec:Beam}.
\end{remark}
	
\subsection{On the normalisation of the eigenvectors}

The values of the coefficients produced by the parametrisation method for the nonlinear mappings and the reduced dynamics directly depend on the selected normalisation of the eigenvectors. As mentioned for instance in \cite{Breden2016}, different normalisations correspond to different radii of convergence in the domain parametrised by the normal coordinates. However, a single validity limit exists in the original physical coordinates. Two equivalent approaches for handling the convergence in both domains are commented on in~\cite{Breden2016}: either setting a given normalisation for the eigenvectors and determining the domain of convergence of the approximation, or selecting an {\em a priori} target domain of convergence and adapting the normalisation of the eigenvectors such that the actual convergence zone corresponds to the prescription. In this contribution, we choose the first approach, as we have deemed it easier to implement for some of the criteria. It should be noted, however, that the image of the manifold that is valid under the parametrisation should be the same, irrespective of which of the two methodologies is adopted. Furthermore, the eigenvectors normalisation is also associated to the decay rate of the parametrisation coefficients, which can impact the numerical conditioning of the method, as pointed out in \cite{Haro,Breden2016}. In general, choosing a scaling that leads to a validity domain of radius close to one seems to ensure good numerical properties to the scheme~\cite{Breden2016}. Nevertheless, different normalisations have been tested in some of the numerical examples, and no numerical problems were found, indicating a certain robustness in this sense for the cases related to nonlinear vibrations.

\subsection{On an {\em a priori} estimation of the maximum forcing amplitude for non-autonomous systems} \label{sec:a_priori_amplitude}

Until this point, all of the presented validity limit criteria are {\em a posteriori} in nature, and, in the case of forced systems, depend on the forcing amplitude for which the parametrisation is computed: they allow for the validation of a certain ROM, computed for a specified load. In daily engineering practice, however, it would be useful to have a criterion allowing for an {\em a priori} estimation of the maximum loading amplitude to which a given system may be submitted. To this purpose, we propose the following approach: first, the unforced and damped system is considered, and its radius of convergence can be calculated following any of the proposed criteria - the undamped system could also be regarded, specially if damping is small and its validity limits have already been computed, {\em e.g.} for a backbone curve calculation. We then suppose that the maximum amplitude for which the normal coordinates defining the parametrisation remain valid is the same for the unforced and the forced systems, and search for a relationship relating the amplitude of the normal coordinates at the peak of the frequency response curves (FRCs) with the loading value in order to deduce the maximum forcing amplitude that the system can be subjected to. Such expressions can be computed using any available analytical approach, such as the parametrisation method itself~\cite{Stabile:morfesym}, or the method of multiple scales~\cite{Nayfeh79}, among many others. Importantly, this relationship shall depend on the type of resonance excited by the forcing (primary or secondary) to be accurate.

To this purpose, expressions for the primary and cubic secondary super- and subharmonic resonances of the Duffing oscillator computed symbolically by the parametrisation method are available in \cite{Stabile:morfesym}, and will be directly used in \cref{sec:DuffingPrimary,sec:superDuffing}. Such solutions, however, are in general not available for high-dimensional systems due to their complexity. We propose then, for a given resonance scenario, to consider a generic form for the reduced dynamics and to employ a semi-analytical approach. For example, considering an arbitrary system in primary resonance reduced to a single master mode using a CNF style, its reduced dynamics is generically written, up to order 3 on the non-autonomous variables and 1 on the autonomous ones, as~\cite{Stabile:morfesym}:
\begin{equation}
	\dot{z}_1 = f_1 z_1 + f_2 z_1^2z_2 + f_3 z_3.
\end{equation}
In this expression, $(z_1,z_2)$ refers to the master mode coordinates, while $(z_3,z_4)$ refers to the forcing~\cite{Vizza:superDPIM,Stabile:morfesym}. One has $z_2=\bar{z}_1$, such that the second equation is the complex conjugate of the first one and is not shown, and $z_3=e^{i\Omega t}$, $z_4=\bar{z}_3$. The coefficients $f_i$ are complex and will be decomposed into real and imaginary parts as $f_i = f_i^R + i f_i^I$. In particular, the first coefficient is always related to the linear characteristics of the system, and is equal to $f_1 = -\xi \omega + i \delta \omega$, with $\omega$ the linear frequency of the master mode, $\xi$ its damping ratio and $\delta=\sqrt{1-\xi^2}$. Additionally, the coefficient $f_3$ is always linear with respect to the forcing amplitude $\kappa$, and will be written as $f_3 = (c_3^R + ic_3^I) \kappa$. Stopping the polynomial order of the expansion at these orders is particularly well-suited to the case of small damping and nonlinearities, where the above reduced dynamics terms dominate. Under these conditions, it is possible to show, as we do in \cref{sec:max_ampl_comp}, that for $c_3^R \neq 0$, the forcing amplitude is related to the parametrisation coefficients by
\begin{equation} \label{eq:kappa_max_gen}
	\kappa = \frac{2 f_1^R \rho + f_2 ^R \frac{\rho^3}{2}}{4c_3^R \sqrt{1 + \alpha^2}},
\end{equation}
with $\alpha = f_3^I/f_3^R$, and a similar equation can be found for the case where $f_3^I \neq 0$, see \cref{sec:max_ampl_comp}. This expression can then be used for any system by inserting the values of the coefficients of the reduced dynamics computed numerically.

\section{Numerical examples} \label{sec:Examples}

We test the criteria on different examples of increasing complexity in order to exemplify how they can be used in practice to determine a safe validity range for the ROMs. For all of the examples a complex normal form style is chosen, unless stated otherwise.

\subsection{The Duffing equation}

In the first example we consider a Duffing oscillator, whose equation of motion is given by
\begin{equation}
	\ddot{u} + 2\xi \omega \dot{u} + \omega^2 u + h u^3 = \kappa \cos{\Omega t}.
\end{equation}
For the numerical examples, $h$ is fixed as $1$ and $\omega$ is taken to be $1.5$, with units compatible.

Since the dimension of the phase space for the Duffing oscillator is two, no reduction actually takes place when the parametrisation method is employed and, if a CNF style of parametrisation is chosen, the technique actually corresponds to the classical normal form approach.

Several scenarios will be considered for the analysis: first, the unforced and undamped case is studied. The main output of the method in this case is the backbone curve of the system; in what follows, two classical resonance scenarios are also studied: primary and superharmonic resonances, in which case damping and forcing are introduced.

\subsubsection{Unforced conservative oscillations} \label{sec:Duffing_consevative}

\begin{figure}[H]
	\centering
	\includegraphics[width=0.5\textwidth]{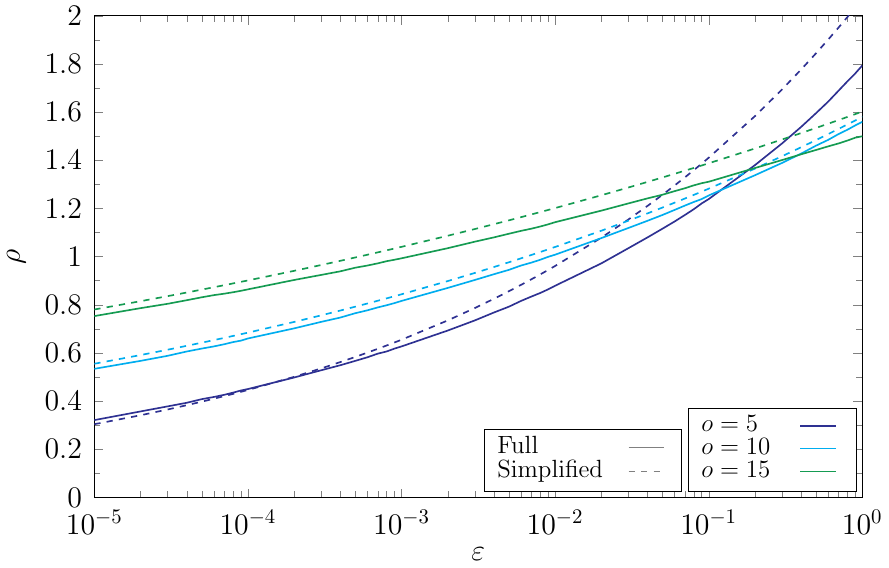}
	\caption{Evolution of the radius of convergence with the tolerance value $\varepsilon$ for the unforced and undamped Duffing oscillator using the invariance equation criterion. The continuous lines indicate results obtained by \cref{eq:inv_eq_error}, while the dashed lines represent the simplification introduced in \cref{eq:inv_eq_error_ANM}. Parameter values are fixed as $\omega = 1.5$ and $h=1$.}
	\label{fig:Duffing cubic conservative unforced - Invariance equation crtiterion}
\end{figure}

First, we consider the invariance equation criterion. In particular, the choice of the error tolerance level plays an important role in this case. \cref{fig:Duffing cubic conservative unforced - Invariance equation crtiterion} presents the evolution of the radius of convergence of the approximation as a function of $\varepsilon$ for different orders of parametrisation. It is interesting to note that the plots get more linear as the parametrisation order increases. For engineering practice, a level of error of 1\%, {\em i.e. $\varepsilon = 0.01$}, is usually acceptable and will be chosen as a reference for further results. The results obtained by applying the simplified formula given by \cref{eq:inv_eq_error_ANM} are also shown, in dashed lines. They present the same qualitative behaviour as the full invariance equation error and serve as a good estimator thereof. 

Continuing, \cref{fig:Duffing cubic conservative unforced - Series criteria} shows the evolution of the radius of convergence through the d'Alembert and Cauchy criteria. In order to apply them, it is necessary to select scalar components of the vectors $\bff(z)$ and $\bfW(z)$ and determine the series radii of convergence for several different values of $\theta$. In particular, the results obtained employing three of them will be depicted: the displacement and velocity, and the first component of the reduced dynamics, denoted $u$, $v$ and $f$ in the figure, respectively. From the plot, it can be seen that the maximum amplitude of the normal coordinates estimated by these three different quantities is virtually the same and that only a small variation in the radius of convergence is observed between the two criteria. It should be noted that for the series criteria, high orders of the parametrisation need to be calculated, such that convergence of the power series is attained. This is not always possible in practice for larger models. Therefore, we also present, in dashed line in \cref{fig:Duffing cubic conservative unforced - Series criteria}, results obtained with parametrisations only up to order 15, obtained by extrapolating the data by means of an exponential fit. A good concordance is obtained between the values predicted by extrapolation and those at order 35.

Finally, the criterion concerned with the singularity of the homological operator is studied. 
The points at which $\det \bfH = 0$ on the complex plane are shown in~\cref{fig:Duffing cubic conservative unforced - Singularity criterion} for three different orders of expansion of the parametrisation. Once these numerical zeros are retrieved, a circle with a radius which is equal to the minimum amplitude among those points is superimposed on the figure for visual guidance. Inside the circle represented with a solid line, the homological operator has no singularities, hence defining the predicted safe amplitude range. For the two presented criteria, convergence with the order of expansion is observed to be monotonic and from above.

\noindent
\begin{minipage}[t]{.48\textwidth}
	\centering
	\includegraphics[width=\textwidth,valign=t]{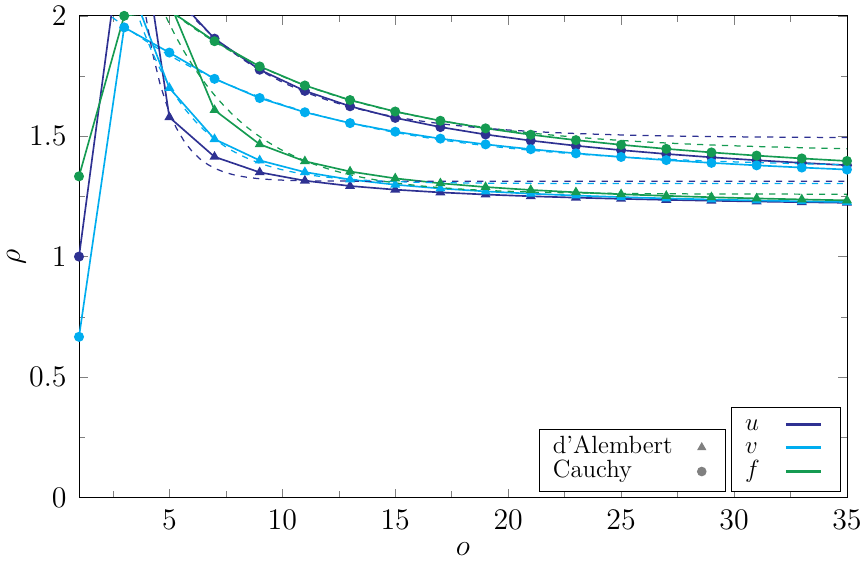}
	\captionof{figure}{Validity limit calculated by the series criteria for the unforced and undamped Duffing oscillator as a function of the parametrisation order. For each of the criteria, the displacement $u$, velocity $v$, and reduced dynamics $f$ power series are studied. Results obtained by extrapolation of parametrisations up to order 15 are also presented, in dashed lines. Parameter values are fixed as $\omega = 1.5$ and $h=1$.}
	\label{fig:Duffing cubic conservative unforced - Series criteria}
\end{minipage}%
\hspace{0.04\textwidth}
\begin{minipage}[t]{.48\textwidth}
	\centering
	\includegraphics[width=\textwidth,valign=t]{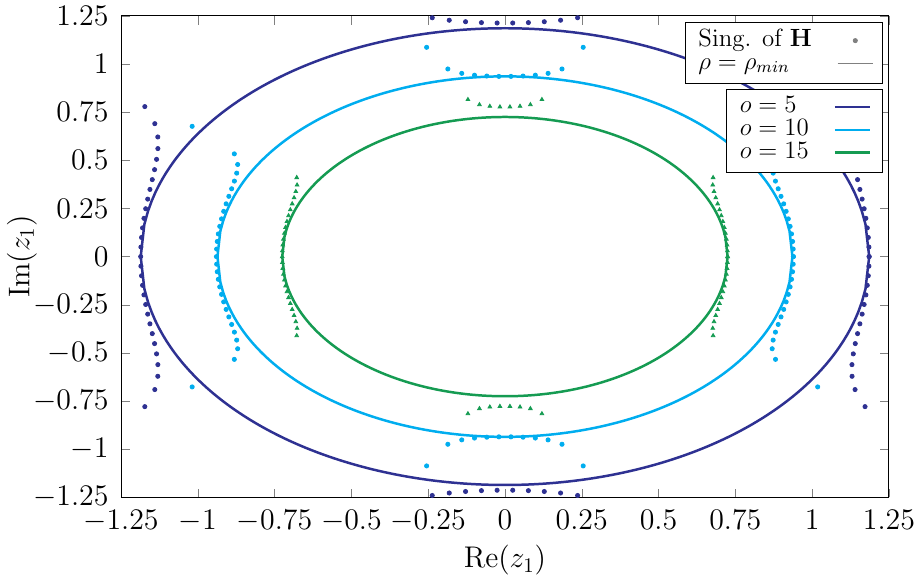}
	\captionof{figure}{Validity limit calculated by the criterion of the singularity of the homological operator for the unforced and undamped Duffing oscillator as a function of the parametrisation order. The lines represent circles of radius equal to the minimum amplitude among points of the same order. Parameter values are fixed as $\omega = 1.5$ and $h=1$.}
	\label{fig:Duffing cubic conservative unforced - Singularity criterion}
\end{minipage}

\begin{figure}[h]
	\centering
	\includegraphics[width=0.5\textwidth]{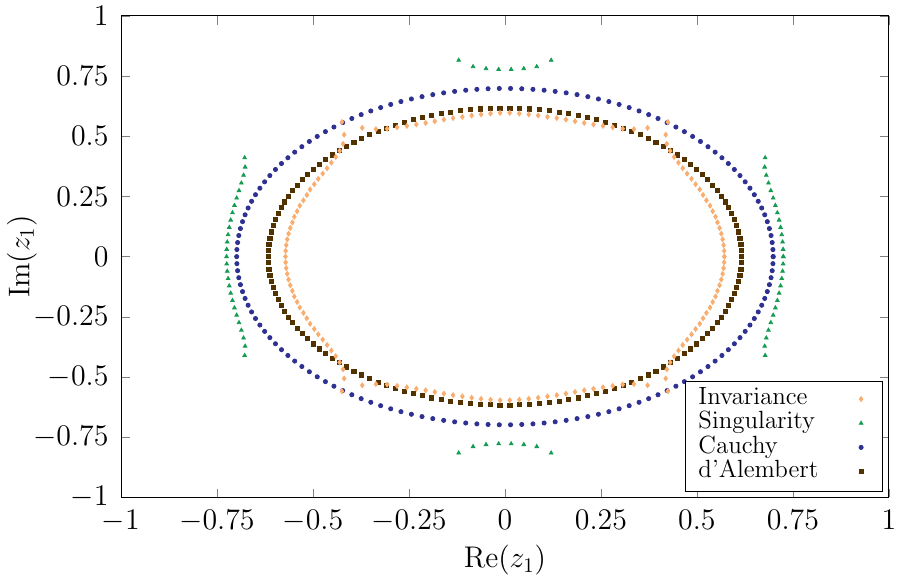}
	\caption{Comparison of the radii of convergence calculated by the tested criteria, for different values of $\theta$. The points correspond to the $z_1$ coordinates in the complex plane for each combination of $\rho$ and $\theta$. For the invariance equation and singularity criteria, the expansion order is chosen as 15, while for the series criteria, it goes up to 35. Parameter values are fixed as $\omega=1.5$ and $h=1$.}
	\label{fig:Duffing cubic conservative unforced - Theta comparison}
\end{figure}

While determining the validity limits, several different values of the phase $\theta$ are tested, and the lowest value of the convergence radius is kept. In order to illustrate the effect of different phases on the computed validity limits, the $z_1$ coordinate on the complex plane for each pair $(\rho,\theta)$ is depicted in \cref{fig:Duffing cubic conservative unforced - Theta comparison}. It is possible to see that, for the series criteria, the phase has little influence on the convergence radius value. For the invariance and singularity criteria, on the other hand, different $\theta$ values have a strong influence on the computed validity limits. This effect is especially pronounced for the singularity criterion, where for a given range of phases, no singularity of the homological operator can be found. Nevertheless, the figure illustrates the converging behaviour of the different criteria in assessing the safe computational region for this particular case.

\begin{figure}[h]
	\centering
	\includegraphics[width=0.39\textwidth]{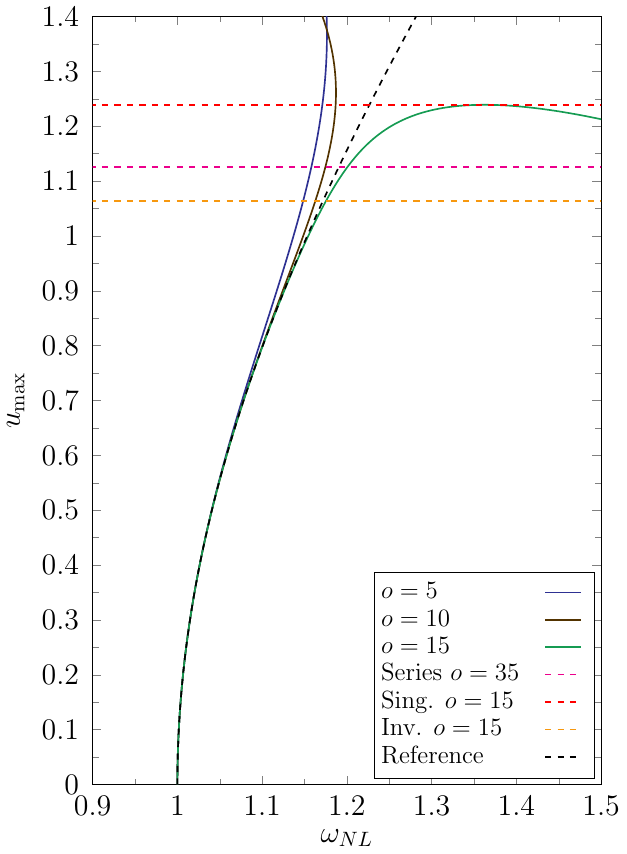}
	\caption{Backbone curves of the conservative unforced Duffing oscillator: CNF with orders $o$=5, 10 and 15 and reference solution~\cite{Salas2014}. Dashed horizontal lines: validity limit criteria obtained by the three methods: series convergence, singularity of the homological operator and error in the invariance equation (Inv.) with $\varepsilon=0.01$. Parameter values are fixed as $\omega = 1.5$ and $h=1$.}
	\label{fig:Duffing cubic conservative unforced - Backbones}
\end{figure}

All the different validity limits encountered are gathered in~\cref{tab:Duffing cubic conservative unforced - Rho values}, both in terms of the amplitude of the normal coordinates, $\rho$, and the maximum physical displacement of the oscillator:

\begin{equation}
	u_{\text{max}} = \max_{0 \leq \theta \leq 2 \pi}{\text{W}_1(\bfz)}.
\end{equation}
Furthermore, the limits are shown graphically together with backbone curves for different orders of parametrisation and a reference analytical solution \cite{Salas2014} in \cref{fig:Duffing cubic conservative unforced - Backbones}. It should be noted that, since the CNF style of parametrisation is employed, the backbone curves can be obtained analytically in a straightforward manner from the reduced dynamics~\cite{BreunungHaller18}. For the series criteria, the value shown corresponds to the minimum of the two in \cref{tab:Duffing cubic conservative unforced - Rho values}. 
The behaviour of the backbone curves approximated by the CNF parametrisation with increasing orders, shows that marginal gains are obtained from $o=15$, which is thus the maximum order represented. Consequently, the validity limits for the singularity and invariance equation criteria are also depicted, both in \cref{tab:Duffing cubic conservative unforced - Rho values} and \cref{fig:Duffing cubic conservative unforced - Backbones} at order 15. This is not the case for the series criteria, whose calculation, by nature, presupposes convergence of the power series representation of the quantities of interest. The radii of convergence are thus given at order 35. However, in order to determine the maximum displacement associated with them, a parametrisation of order 15 is chosen. This is done because, due to the high order of expansion, slight variations in the normal variables amplitude can correspond to large variations of the physical vibrations associated to them. Also, this approach is in line with the idea of applying the extrapolation procedure for high dimensional systems, since in this scenario one does not have access to the high order parametrisation and needs to employ the highest order available. This will be the case for all of the examples presented in the text.

\begin{table}[h]
	\npdecimalsign{.}
	\nprounddigits{3}
	\centering
	\begin{tabular}{c|c|c|c|c|c}
		Criterion & Cauchy & \text{d'Alembert} & Singularity & Invariance & Invariance simpl. \\ \hline
		$\rho$ & 1.361 & 1.224 & 1.450 & 1.143 & 1.202 \\
		$u_{\text{max}}$ & 1.215 & 1.126 & 1.240 & 1.063 & 1.109 \\
	\end{tabular}
	\caption{Validity limits for the conservative unforced Duffing oscillator. Values are given in both the amplitude of the normal coordinates and maximum physical displacement of the oscillator. Orders of expansion are 15 for the singularity and invariance equation criteria and 35 for the series criteria. The tolerance for the invariance equation is chosen as $\varepsilon = 1\%$. Parameter values are fixed as $\omega = 1.5$ and $h=1$.}
	\label{tab:Duffing cubic conservative unforced - Rho values}
\end{table}

As can be seen in \cref{fig:Duffing cubic conservative unforced - Backbones}, the calculated validity limits, although presenting different values, intuitively give a good indication of the zone at which the backbone curves obtained by the parametrisation method, deviate from the reference solution. Furthermore, by comparing the values of the convergence radii given in \cref{fig:Duffing cubic conservative unforced - Series criteria,fig:Duffing cubic conservative unforced - Singularity criterion,fig:Duffing cubic conservative unforced - Invariance equation crtiterion} for different orders, it is possible to note that the series and singularity criteria converge from above, while the invariance equation one converges from below. As such, their values can be interpreted as bounds to the region where the parametrisation loses convergence. In practical engineering analysis, this would be a zone in which results should be taken with caution.

\subsubsection{Primary resonance} \label{sec:DuffingPrimary}

For this case, the damping ratio is fixed at $\xi = 0.02$. As proposed, for example, in~\cite{opreni22high,Vizza:superDPIM}, the ROM for the forced-damped case is derived for a single value of the excitation frequency, and then used to compute the FRC in the vicinity of the selected resonance case. The selected excitation frequency retained to compute the ROM is denoted as the expansion point, which is here set at the undamped natural frequency of the oscillator: $\Omega = \omega$. Additionally, three cases regarding forcing amplitude are considered: $\kappa = 0.1$, $\kappa = 0.175$, and $\kappa = 0.25$. A detailed study of each criterion will be made for $\kappa = 0.175$, while only the final results together with the FRCs at different orders will be given for the other two forcing values.

\begin{figure}[h]
	\centering
	\includegraphics[width=0.5\textwidth]{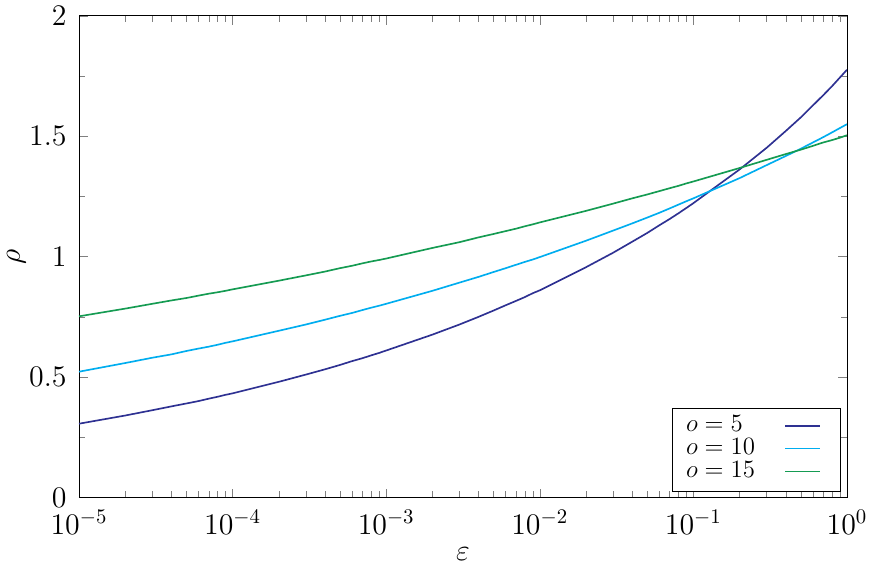}
	\caption{Evolution of the radius of convergence with the tolerance value for the damped forced Duffing oscillator in primary resonance using the invariance equation criterion. Parameter values are fixed as $\omega = 1.5$, $h=1$, $\xi = 0.02$ and $\kappa = 0.175$.}
	\label{fig:Duffing cubic damped primary - kappa = 0.175 - Invariance equation criterion}
\end{figure}

Since the forcing is present in this scenario, several values of the phase $\phi$ should also be considered when sampling values of the normal coordinates. Specifically, the values $\phi \in \{0, \frac{\pi}{2}, \pi, \frac{3\pi}{2}\}$ are chosen. \cref{fig:Duffing cubic damped primary - kappa = 0.175 - Invariance equation criterion} shows the evolution of the convergence radius obtained by the invariance equation criterion with the imposed tolerance value $\varepsilon$. Once again, as the order of parametrisation increases, the curves tend to be more linear.

\begin{figure}[h]
\begin{minipage}[t]{.48\textwidth}
	\centering
	\includegraphics[width=\textwidth,valign=t]{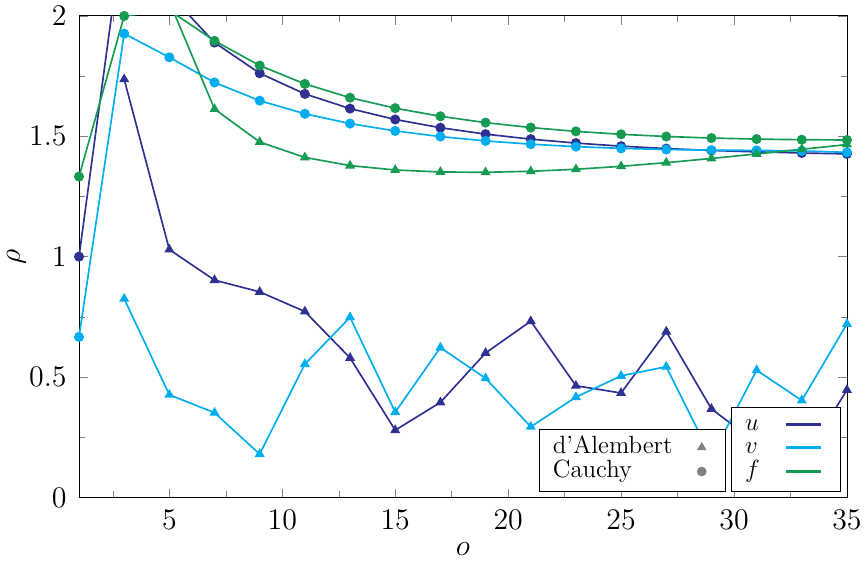}
	\captionof{figure}{Validity limit calculated by the series criteria for the damped Duffing oscillator forced in primary resonance as a function of the parametrisation order $o$. For each of the criteria, the displacement, velocity and reduced dynamics power series are studied. Parameter values are set as $\omega = 1.5$, $h=1$, $\xi = 0.02$, $\kappa = 0.175$ and $\phi = \frac{\pi}{2}$.}
	\label{fig:Duffing cubic damped primary - kappa = 0.175 - Series criteria}
\end{minipage}%
\hspace{0.04\textwidth}
\begin{minipage}[t]{.48\textwidth}
	\centering
	\includegraphics[width=\textwidth,valign=t]{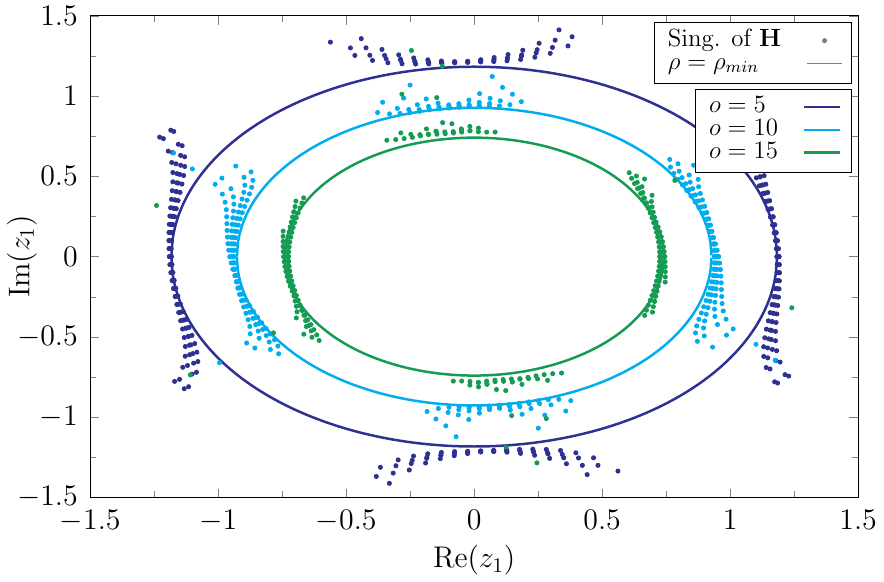}
	\captionof{figure}{Validity limit calculated by the criterion of the singularity of the homological operator for the damped Duffing oscillator in primary resonance as a function of the parametrisation order. The lines represent circles of radius equal to the minimum amplitude among points of the same order. Parameter values are fixed as $\omega = 1.5$, $h=1$, $\xi = 0.02$ and $\kappa = 0.175$.}
	\label{fig:Duffing cubic damped primary - kappa = 0.175 - Singularity crtiterion}
\end{minipage}
\end{figure}

The evolution of the radii of convergence given by the series criteria is depicted in \cref{fig:Duffing cubic damped primary - kappa = 0.175 - Series criteria}. Note that four different values of the phase $\phi$ are used in the computation ($\phi=0,\pi/2,\pi, 3\pi/2$) of the validity limit associated with this criterion, but only the case showing $\phi = \pi/2$, corresponding to the minimal final value of the validity limit, is displayed in the figure. In this case, d'Alembert's criterion fails to yield suitable results: the series for displacement and velocity do not converge, whereas the series for the reduced dynamics exhibits an increasing radius of convergence at higher polynomial orders. For the Cauchy rule, the three series yield nearly identical values, with monotonic convergence.

Finally, \cref{fig:Duffing cubic damped primary - kappa = 0.175 - Singularity crtiterion} shows the results for the singularity criterion. Again, the four values of the phase $\phi$ are used and now displayed, hence giving rise to four sets of vanishing points per order. The circle, used for visual guidance, takes the minimum of the four cases.
Qualitatively, the observed evolution is similar to that of the unforced case, with the points at which the homological operator becomes singular being almost the same for the different tested values of $\phi$.

\begin{figure}[h]
	\centering
	\begin{subfigure}{0.65\textwidth}
		\centering
		\includegraphics[width=\textwidth]{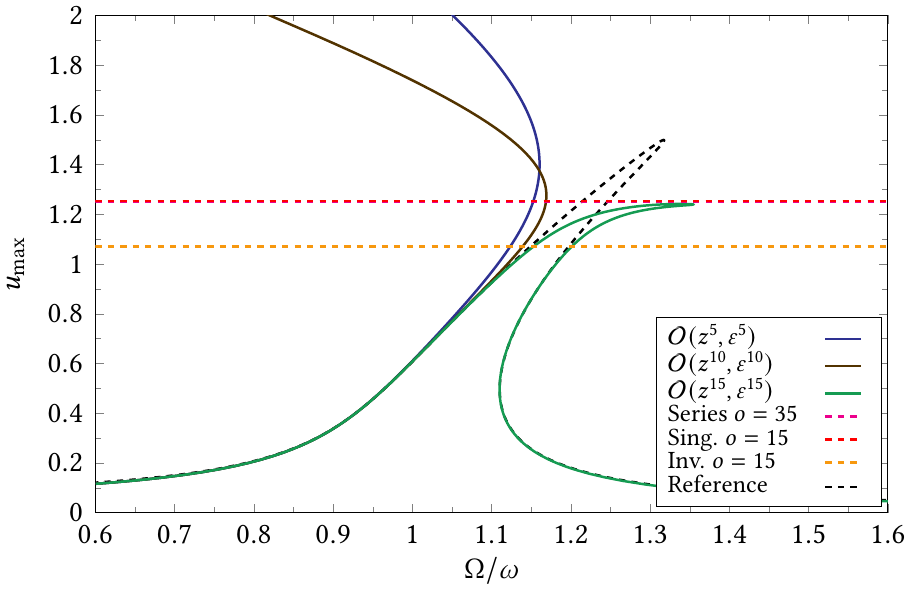}
		\caption{$\kappa = 0.175$}
		\label{fig:Duffing cubic damped primary - kappa = 0.175 - FRCs}
	\end{subfigure}
	
	\begin{subfigure}{0.48\textwidth}
		\centering
		\includegraphics[width=\textwidth]{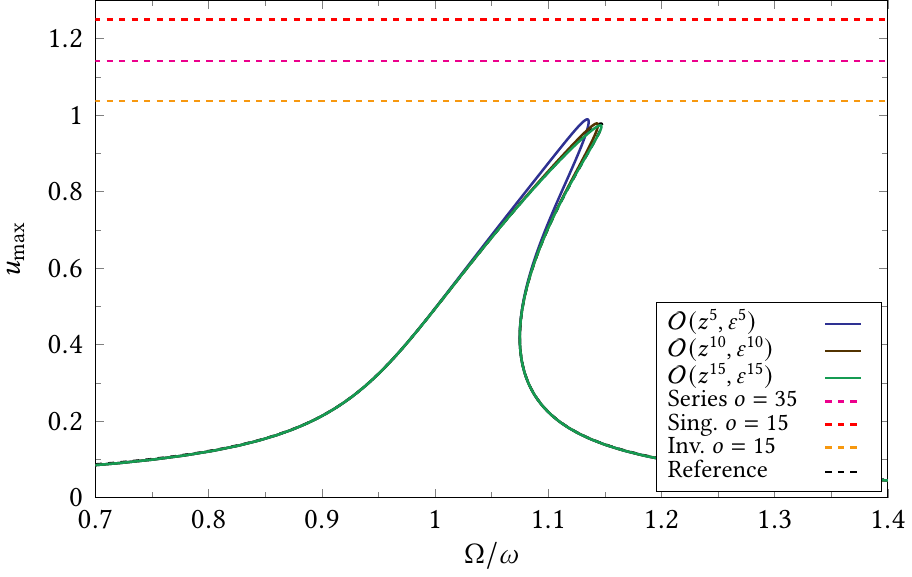}
		\caption{$\kappa = 0.1$}
		\label{fig:Duffing cubic damped primary - kappa = 0.1 - FRCs}
	\end{subfigure}
	\begin{subfigure}{0.48\textwidth}
		\centering
		\includegraphics[width=\textwidth]{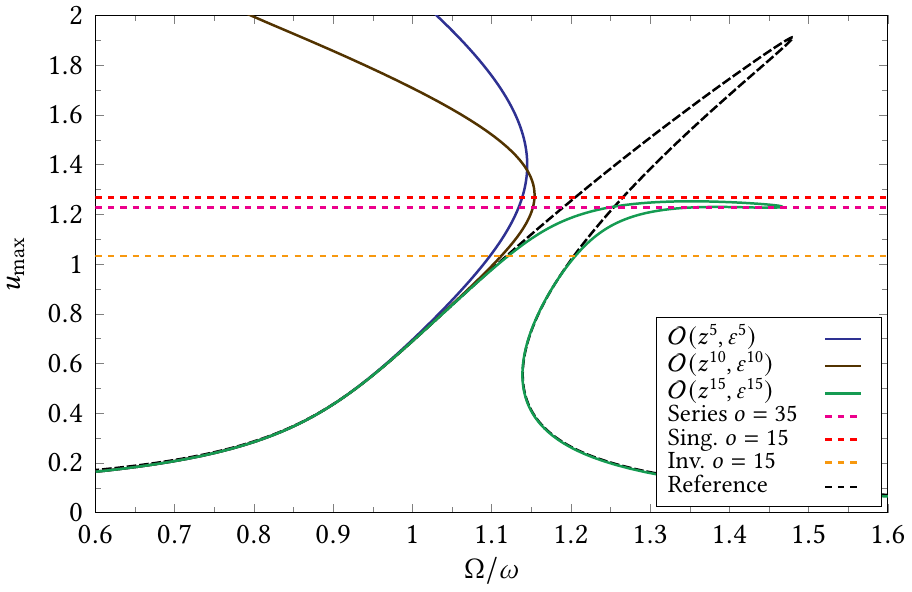}
		\caption{$\kappa = 0.25$}
		\label{fig:Duffing cubic damped primary - kappa = 0.25 - FRCs}
	\end{subfigure}
	\caption{Frequency response curves for the primary resonance of the Duffing oscillator. Different orders of parametrisation are compared with a reference solution computed by means of numerical continuation, where the matcont package~\cite{dhooge2004matcont}, has been used. The validity limits calculated by the different proposed convergence criteria are also shown. Parameter values are fixed as $\omega = 1.5$, $h=1$ and $\xi = 0.02$.}
	\label{fig:Duffing cubic damped primary - FRCs}
\end{figure}

The amplitudes of the normal variables obtained by each criterion, as well as the associated maximum physical displacements, are given in~\cref{tab:Duffing cubic damped primary - Rho values} for $\kappa = 0.175$. Additionally, the FRCs for all values of forcing amplitude are depicted in~\cref{fig:Duffing cubic damped primary - FRCs}. For $\kappa = 0.1$, displayed in~\cref{fig:Duffing cubic damped primary - kappa = 0.1 - FRCs}, all the computed validity limits are well above the peak of the FRCs, ensuring a perfect convergence of the solutions produced by the ROM for that forcing amplitude. For the other two forcing scenarios, the convergence radii show a good agreement with the region at which the ROMs fails to converge. Once again, with the selected tolerance value for the invariance equation criterion ($\varepsilon=0.01$), it seems to yield a lower bound to the zone at which ROM results deviate from those of the FOM. 

\begin{table}[h]
	\npdecimalsign{.}
	\nprounddigits{3}
	\centering
	\begin{tabular}{c|c|c|c|c}
		Criterion & Cauchy & \text{d'Alembert} & Singularity & Invariance \\ \hline
		$\rho$ & 1.427 & 1.465 & 1.482 & 1.142 \\
		$u_{\text{max}}$ & 1.251 & 1.257 & 1.257 & 1.071\\
	\end{tabular}
	\caption{Validity limits for the damped Duffing oscillator in primary resonance condition. Values are given in both the amplitude of the normal coordinates $\rho$ and maximum physical displacement $u_{\text{max}}$ of the oscillator. Orders of expansion are 15 for the singularity and invariance equation criteria and 35 for the series criteria. The tolerance for the invariance equation is chosen as $\varepsilon = 1\%$. Parameter values are fixed as $\omega = 1.5$, $h=1$, $\xi = 0.02$ and $\kappa = 0.175$.}
	\label{tab:Duffing cubic damped primary - Rho values}
\end{table}

To conclude with this example, the approach proposed in~\cref{sec:a_priori_amplitude} is illustrated for this particular case of the primary resonance of the Duffing oscillator. By applying the invariance equation criterion for the unforced and damped Duffing oscillator, the radius of convergence is found to be $\rho = 1.144$. For this case, the parametrisation coefficients can be calculated analytically, see \cite{Stabile:morfesym}, and by replacing them into \cref{eq:max_kappa_fI}, one finds:
\begin{equation}
	\kappa = 2\omega^2\xi\rho.
\end{equation}
The intuitively expected behaviour is observed: as damping increases, the maximum load the system can sustain also increases, in a linear manner. For the Duffing oscillator in particular, due to the simplicity of the system, it is also possible to consider a form of the reduced dynamics which is more general than that used to derive \cref{eq:max_kappa_fI}, while still yielding an analytical expression relating $\rho$ and $\kappa$ at the peak. This formula, not detailed here for the sake of brevity, can be found in~\cite{Stabile:morfesym}, and the more accurate relationship then reads
\begin{equation}
	\kappa = \frac{32 \omega^4 \xi \rho}{3 h \rho^2 - 16 \omega^2}.
\end{equation}
The first of these equations yields $\kappa = 0.103$, while the second gives $\kappa  = 0.116$. We see that the values are close, and thus, indeed, the form of the reduced dynamics chosen in \cref{sec:a_priori_amplitude} is a good approximation in this example. By inspection of \cref{fig:Duffing cubic damped primary - kappa = 0.1 - FRCs,fig:Duffing cubic damped primary - kappa = 0.175 - FRCs}, it becomes clear that these values are a good estimation of the maximum loading amplitude that the system can be subjected to. Hence, the validity limit found in the unforced case can be easily used to deduce the maximum forcing amplitude for which convergence is achieved.

\subsubsection{Superharmonic resonance} \label{sec:superDuffing}

For this case, the damping ratio is fixed at $\xi = 0.002$, and the expansion forcing frequency is chosen to coincide with one third of the undamped natural frequency of the oscillator, $\Omega = \frac{\omega}{3}$. Once again, three cases regarding forcing amplitude are considered: $\kappa = 0.5$, $\kappa = 0.75$ and $\kappa = 1$, in order to illustrate what happens for moderate, high and extremely high forcing.

\begin{figure}[h]
	\centering
	\includegraphics[width=0.5\textwidth]{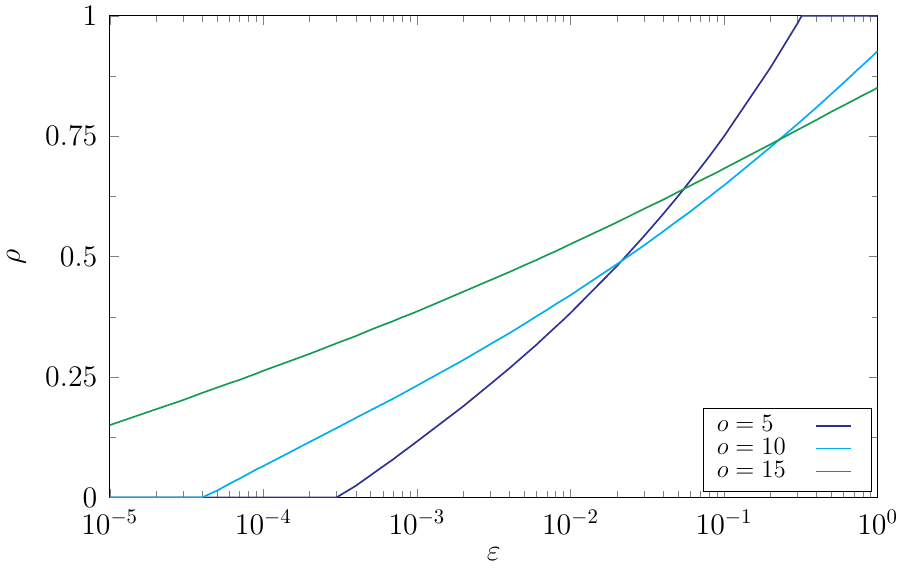}
	\caption{Evolution of the radius of convergence with the tolerance value for the damped Duffing oscillator in superharmonic resonance using the invariance equation criterion. Parameter values are fixed as $\omega = 1.5$, $h=1$, $\xi = 0.002$ and $\kappa = 0.75$.}
	\label{fig:Duffing cubic damped superharmonic - kappa = 0.75 - Invariance equation criterion}
\end{figure}

\begin{figure}[h]
\begin{minipage}[t]{.48\textwidth}
	\centering
	\includegraphics[width=\textwidth,valign=t]{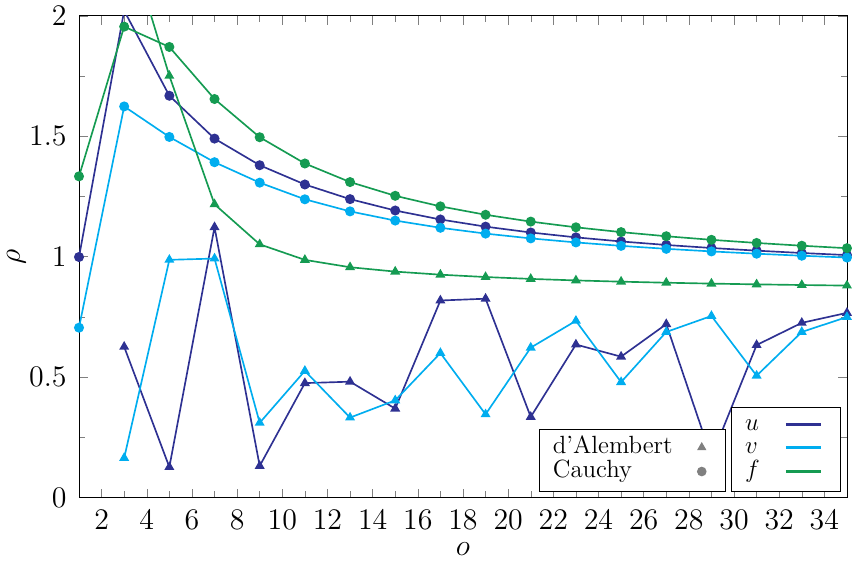}
	\captionof{figure}{Validity limit calculated by the series criteria for the damped Duffing oscillator in superharmonic resonance as a function of the parametrisation order. For each of the criteria, the displacement, velocity and reduced dynamics power series are studied. Parameter values are fixed as $\omega = 1.5$, $h=1$, $\xi = 0.002$, $\kappa = 0.75$ and $\phi=\frac{\pi}{2}$.}
	\label{fig:Duffing cubic damped superharmonic - kappa = 0.75 - Series criteria}
\end{minipage}%
\hspace{0.04\textwidth}
\begin{minipage}[t]{.48\textwidth}
	\centering
	\includegraphics[width=\textwidth,valign=t]{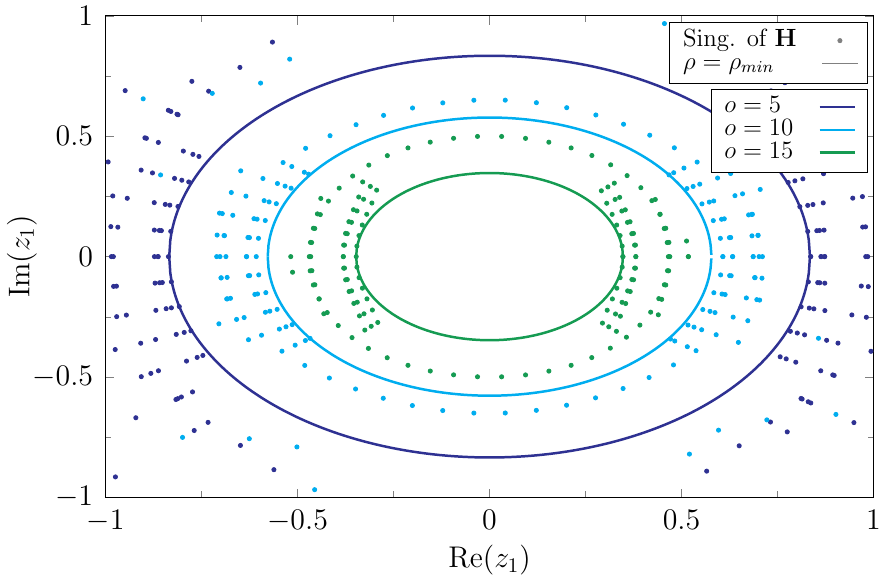}
	\captionof{figure}{Validity limit calculated by the criterion of the singularity of the homological operator for the damped Duffing oscillator in superharmonic resonance as a function of the parametrisation order. The lines represent circles of radius equal to the minimum amplitude among points of the same order. Parameter values are fixed as $\omega = 1.5$, $h=1$, $\xi = 0.002$ and $\kappa = 0.75$.}
	\label{fig:Duffing cubic damped superharmonic - kappa = 0.75 - Singularity crtiterion}
\end{minipage}
\end{figure}

The results for the invariance equation, series and singularity criteria are presented in \cref{fig:Duffing cubic damped superharmonic - kappa = 0.75 - Invariance equation criterion,fig:Duffing cubic damped superharmonic - kappa = 0.75 - Series criteria,fig:Duffing cubic damped superharmonic - kappa = 0.75 - Singularity crtiterion}. By inspecting \cref{fig:Duffing cubic damped superharmonic - kappa = 0.75 - Invariance equation criterion}, one can remark that the choice of the error tolerance value $\varepsilon$ has a more important effect than in the previous cases, as evidenced by the high slope of the curves in the figure. It should be noted, again, that in \cref{fig:Duffing cubic damped superharmonic - kappa = 0.75 - Series criteria} only the data series corresponding to $\phi=\pi/2$ is displayed, although four values of $\phi$ are used in the computation of the validity limit. Also, for the singularity criterion, where vanishing points are represented for four different phases $\phi$, the influence of the phase difference between the response and the forcing is more significant than in the primary resonance case.

\begin{table}[h]
	\npdecimalsign{.}
	\nprounddigits{3}
	\centering
	\begin{tabular}{c|c|c|c|c}
		Criterion & Cauchy & \text{d'Alembert} & Singularity & Invariance \\ \hline
		$\rho$ & 0.996 & 0.880 & 0.762 & 0.526 \\
		$u_{\text{max}}$ & 1.009 & 0.904 & 0.830 & 0.602\\
	\end{tabular}
	\caption{Validity limits for the damped Duffing oscillator forced at the superharmonic resonance. Values are given in both the amplitude of the normal coordinates $\rho$ and maximum physical displacement $u_{\text{max}}$ of the oscillator. Orders of expansion are 15 for the singularity and invariance equation criteria and 35 for the series criteria. The tolerance for the invariance equation is chosen as $\varepsilon = 1\%$. Parameter values are fixed as $\omega = 1.5$, $h=1$, $\xi = 0.002$ and $\kappa = 0.75$.}
	\label{tab:Duffing cubic damped superharmonic - Rho values}
\end{table}

\begin{figure}[h]
	\centering
	\begin{subfigure}{0.65\textwidth}
		\centering
		\includegraphics[width=\textwidth]{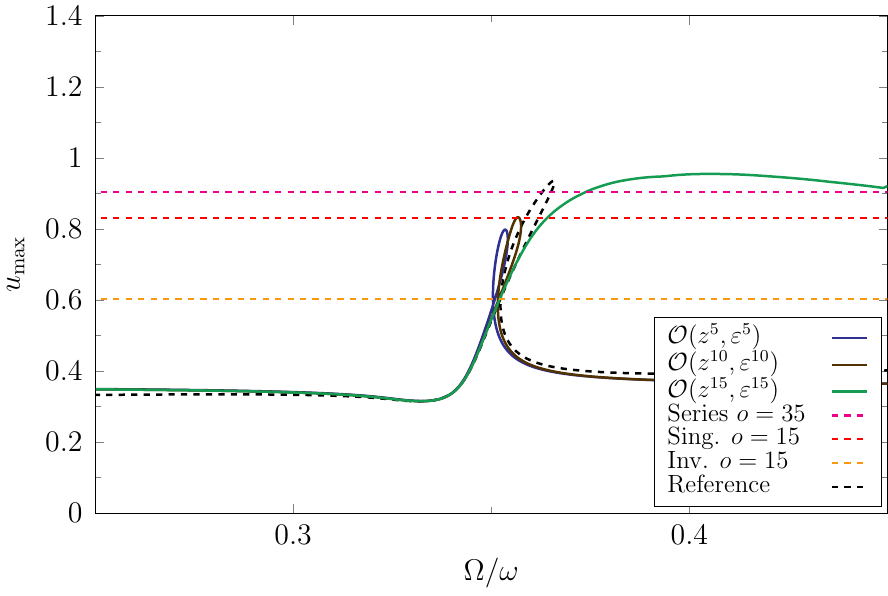}
		\caption{$\kappa = 0.75$}
		\label{fig:Duffing cubic damped superharmonic - kappa = 0.75 - FRCs}
	\end{subfigure}
	
	\begin{subfigure}{0.48\textwidth}
		\centering
		\includegraphics[width=\textwidth]{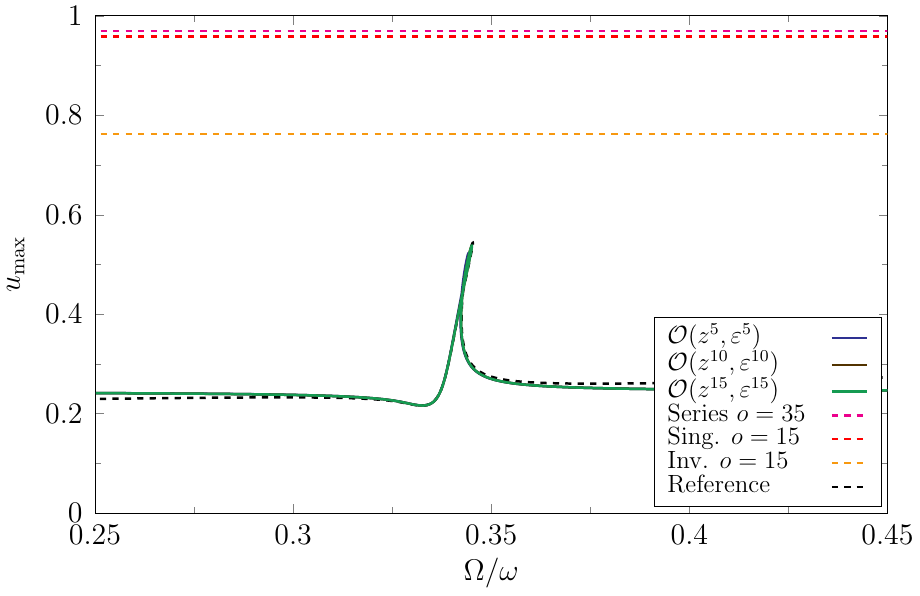}
		\caption{$\kappa = 0.5$}
		\label{fig:Duffing cubic damped superharmonic - kappa = 0.5 - FRCs}
	\end{subfigure}
	\begin{subfigure}{0.48\textwidth}
		\centering
		\includegraphics[width=\textwidth]{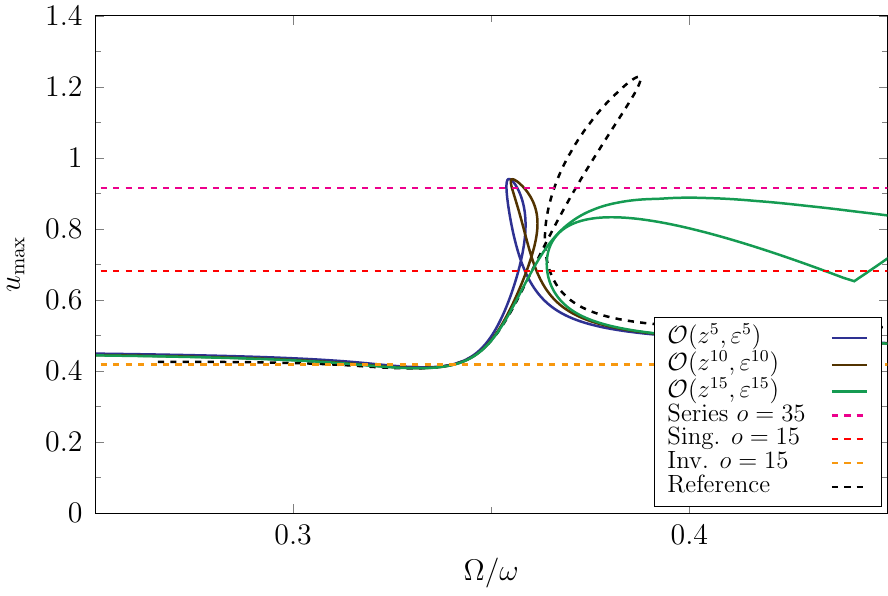}
		\caption{$\kappa = 1$}
		\label{fig:Duffing cubic damped superharmonic - kappa = 1 - FRCs}
	\end{subfigure}
	\caption{Frequency response curves for the superharmonic resonance of the Duffing oscillator. Different orders of parametrisation are compared with a reference solution computed by numerical continuation employing package matcont \cite{dhooge2004matcont}. The validity limits calculated by the different proposed convergence criteria are also shown. Parameter values are fixed as $\omega = 1.5$, $h=1$ and $\xi = 0.002$.}
	\label{fig:Duffing cubic damped superharmonic - FRCs}
\end{figure}

The values of the validity limits obtained by the various criteria for $\kappa = 0.75$ are detailed in \cref{tab:Duffing cubic damped superharmonic - Rho values}, and compared to the FRCs for the three forcing values in~\cref{fig:Duffing cubic damped superharmonic - FRCs}. In this example, the results indicate perfect convergence for the low forcing value, and the intuitively expected convergence region is retrieved for $\kappa = 0.75$. One can observe that, for the higher forcing with $\kappa=1$, the invariance criterion indicates amplitudes for the loss of convergence that encompass points on the FRC for which convergence would be expected due to the agreement with the continuation solution. This can indicate that the reduced-order model correctly captures steady-state oscillatory behaviour associated with the FRC, but is inadequate for other parts of the manifold, which do not correspond to periodic orbits and thus are not represented on the FRC. This result also underlines that the choice of $\varepsilon$ is of importance, and that setting $\varepsilon=0.01$ generally provides a lower bound, which must be complemented with an upper bound, here provided by the two other criteria.

For this resonance scenario, we also estimate the maximum forcing value by using the computed validity limit for the damped backbone curve. The expression relating the forcing amplitude and the peak amplitude of the FRC normal coordinates, in the case of the superharmonic forcing, can be easily derived from the analytical results provided in~\cite{Stabile:morfesym}. Up to the third order, it yields
\begin{equation}
	\kappa = \frac{4}{9} \sqrt[3]{\frac{\omega^8 \left( 9 \xi^2 + 16 \right)^{\nicefrac{3}{2}} \delta \xi \rho}{h}}.
\end{equation}
This relationship gives $\kappa = 0.691$ for $\rho = 1.144$, which by inspection of \cref{fig:Duffing cubic damped superharmonic - kappa = 0.5 - FRCs,fig:Duffing cubic damped superharmonic - kappa = 0.75 - FRCs} is a good estimation of the loading limit for this case.

\subsection{A 2-DOF system} \label{sec:2DOF}

In this section, we consider a 2-DOF system having as equations of motion:
\begin{align}
	\ddot{u}_1 + 2 \xi_1 \omega_1 \dot{u}_1 + \omega_1^2 u_1 + h^1_{111} u_1^3 + 3h^2_{111} u_1^2 u_2 &= \kappa \cos{\Omega t}, \\
	\ddot{u}_2 +  2 \xi_2 \omega_2 \dot{u}_2 + \omega_2^2 u_2 + h^2_{111} u_1^3 &= 0.
\end{align}
This example is chosen in order to test the proposed validity criteria for a case where the number of degrees of freedom in the model is larger than the number of normal coordinates, but that is still of small size. For the numerical examples, $h^1_{111}$, $h^2_{111}$ and $\omega_1$ are set to~$1$. 
The unforced case is first investigated with $\kappa = \xi_1 = \xi_2 = 0$. Two values for $\omega_2$ are chosen: $\omega_2 = 1.57$, causing the behaviour of the system to initially be of the hardening type, and then switch to the softening, and $\omega_2 = 0.637$, causing hardening for all displacement amplitudes.

\begin{figure}[h]
	\centering
	\begin{subfigure}{0.49\textwidth}
		\includegraphics[width=\textwidth]{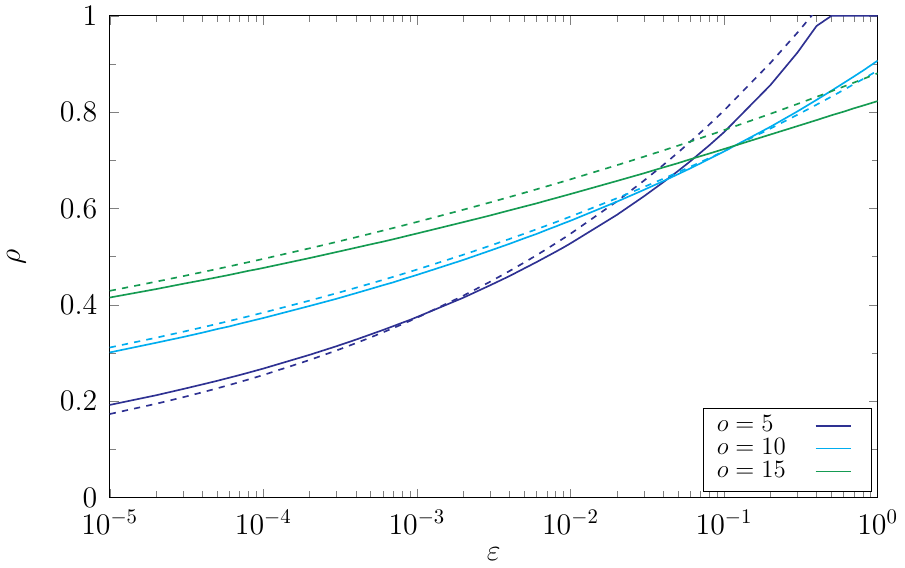}
		\caption{$\omega_2 = 1.57$}
		\label{fig:2 DOF cubic conservative unforced - omega2 = 1.57 - Invariance equation criterion}
	\end{subfigure}
	\begin{subfigure}{0.49\textwidth}
		\includegraphics[width=\textwidth]{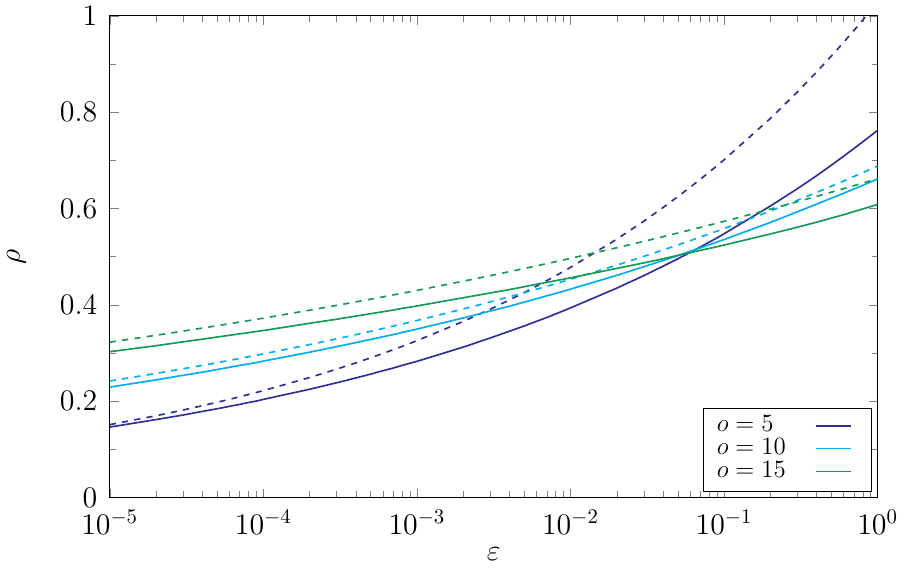}
		\caption{$\omega_2 = 0.637$}
		\label{fig:2 DOF cubic conservative unforced - omega2 = 0.637 - Invariance equation criterion}
	\end{subfigure}
	\caption{Evolution of the radius of convergence with the tolerance value for the unforced and undamped 2-DOF oscillator using the invariance equation criterion. The continuous lines indicate results obtained by \cref{eq:inv_eq_error}, while the dashed ones represent the simplification introduced in \cref{eq:inv_eq_error_ANM}. Parameter values are fixed as $\omega_1 = h^1_{111} = h^2_{111} = 1$.}
	\label{fig:2 DOF cubic conservative unforced - Invariance equation criterion}
\end{figure}

The tolerance versus radius of convergence curves obtained for the invariance equation criterion are depicted in \cref{fig:2 DOF cubic conservative unforced - Invariance equation criterion}. An excellent agreement is observed between the simplified and full ways of evaluating the invariance equation, especially when higher orders of asymptotic developments are employed. From the figure, an additional indicator can be identified for assessing the necessity of a higher-order parametrisation: a marked deviation between the full invariance equation and its simplified counterpart signals the need to compute a higher-order expansion.

\begin{figure}[h]
	\centering
	\begin{subfigure}{0.49\textwidth}
		\includegraphics[width=\textwidth]{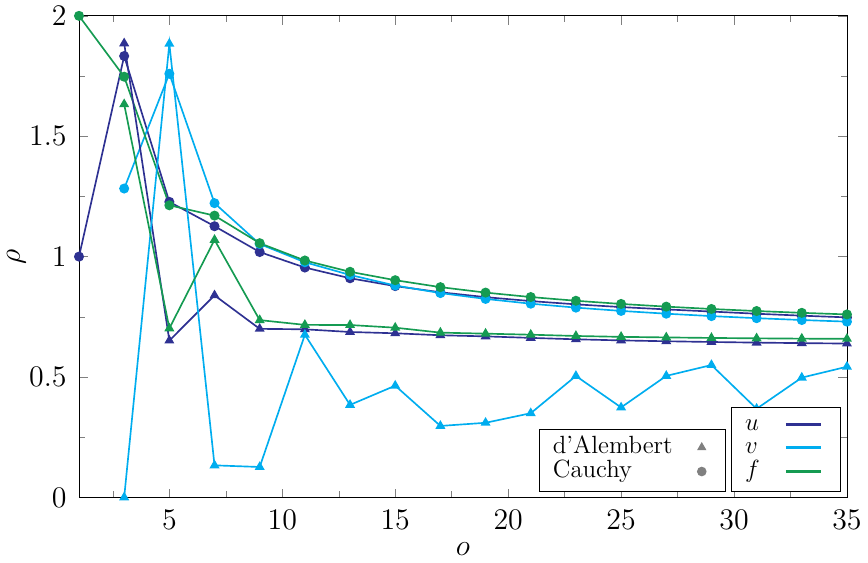}
		\caption{$\omega_2 = 1.57$}
		\label{fig:2 DOF cubic conservative unforced - omega2 = 1.57 - Series criteria}
	\end{subfigure}
	\begin{subfigure}{0.49\textwidth}
		\includegraphics[width=\textwidth]{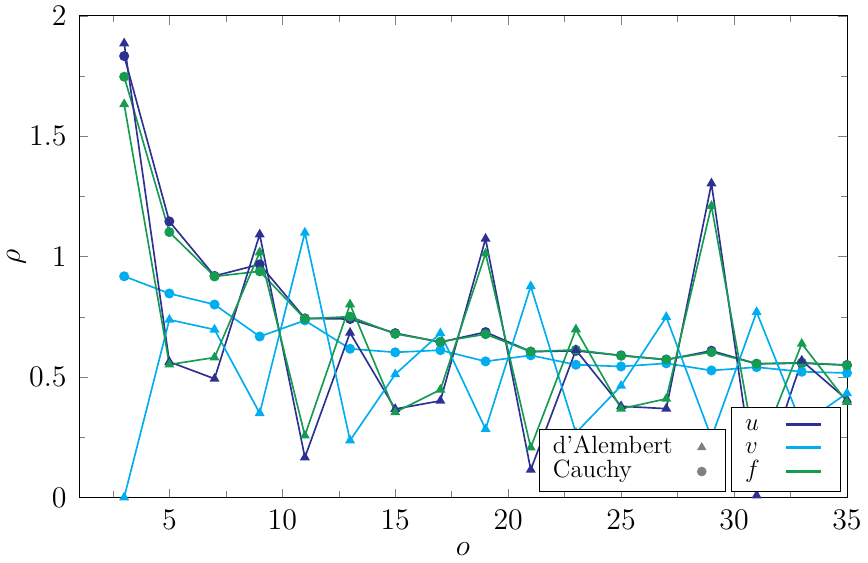}
		\caption{$\omega_2 = 0.637$}
		\label{fig:2 DOF cubic conservative unforced - omega2 = 0.637 - Series criteria}
	\end{subfigure}
	\caption{Validity limit calculated by the series criteria for the unforced and undamped 2-DOF oscillator as a function of the parametrisation order. For each of the criteria, the displacement, velocity and reduced dynamics power series are studied. Parameter values are fixed as $\omega_1 = h^1_{111} = h^2_{111} = 1.0$.}
	\label{fig:2 DOF cubic conservative unforced - Series criteria}
\end{figure}

\begin{figure}[h]
	\centering
	\begin{subfigure}{0.49\textwidth}
		\includegraphics[width=\textwidth]{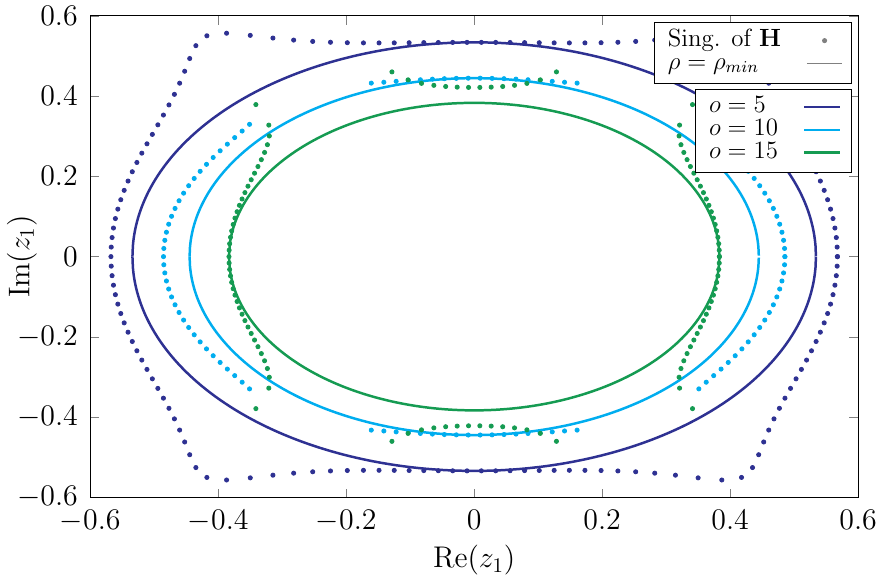}
		\caption{$\omega_2 = 1.57$}
		\label{fig:2 DOF cubic conservative unforced - omega2 = 1.57 - Singularity criterion}
	\end{subfigure}
	\begin{subfigure}{0.49\textwidth}
		\includegraphics[width=\textwidth]{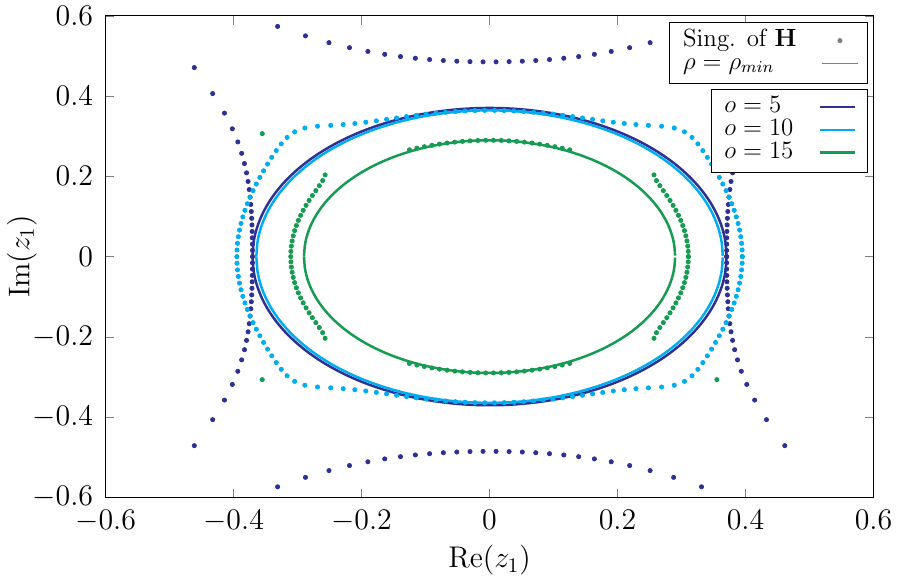}
		\caption{$\omega_2 = 0.637$}
		\label{fig:2 DOF cubic conservative unforced - omega2 = 0.637 - Singularity criterion}
	\end{subfigure}
	\caption{Validity limit calculated by the singularity of the homological operator for the unforced and undamped 2-DOF oscillator as a function of the parametrisation order. The lines represent circles of radius equal to the minimum amplitude among points of the same order. Parameter values are fixed as $\omega_1 = h^1_{111} = h^2_{111} = 1.0$.}
	\label{fig:2 DOF cubic conservative unforced - Singularity criterion}
\end{figure}

The results for the series criteria are presented in \cref{fig:2 DOF cubic conservative unforced - Series criteria} for both values of the second oscillator's linear frequency. In this case, the series $u$ and $v$ are associated to the oscillator's $u_1$ and $v_1$ DOFs. For the case $\omega_2 = 0.637$, the d'Alembert criterion completely fails to converge. The Cauchy rule, on the other hand, seems to converge to similar values for all series, but monotonicity is lost. The application of the singularity criterion in this example differs from the other cases. Since the homological operator is represented by a rectangular matrix, a subset of its rows must be selected in order to compute its determinant. In particular, the rows corresponding to the $u_1$ degree of freedom have been selected, as being related to the master coordinate. The results are presented in \cref{fig:2 DOF cubic conservative unforced - Singularity criterion}. Once again, an interesting result can be observed for the case $\omega_2 = 0.637$, depicted in \cref{fig:2 DOF cubic conservative unforced - omega2 = 0.637 - Singularity criterion}: between orders 5 and 10 there is almost no change in the radius of convergence predicted by the criterion. However, the expected decreasing behaviour is found when going up to order 15. 

\begin{figure}[h]
	\centering
	\begin{subfigure}{0.39\textwidth}
		\includegraphics[width=\textwidth]{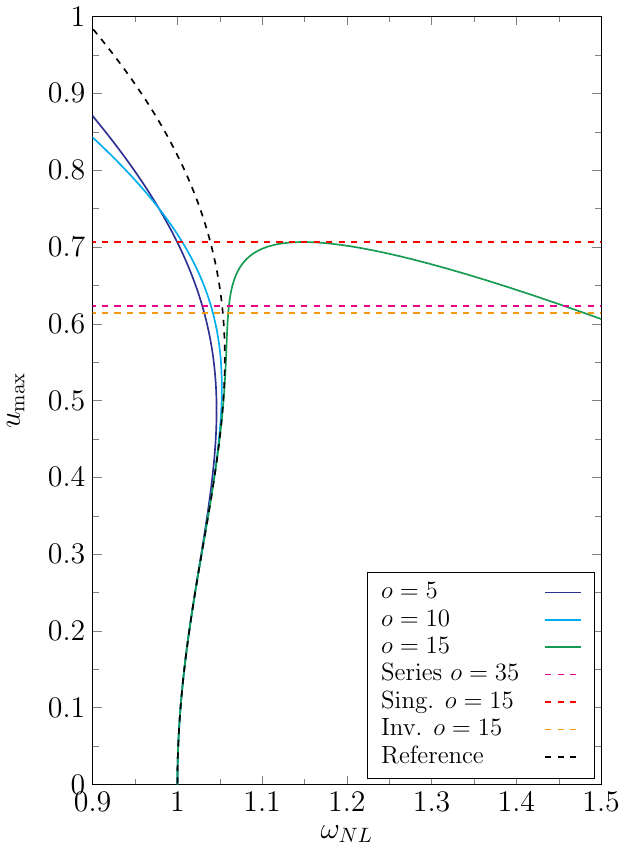}
		\caption{$\omega_2 = 1.57$}
		\label{fig:2 DOF cubic conservative unforced - omega2 = 1.57 - Backbones}
	\end{subfigure}
	\hspace{0.05\textwidth}
	\begin{subfigure}{0.39\textwidth}
		\includegraphics[width=\textwidth]{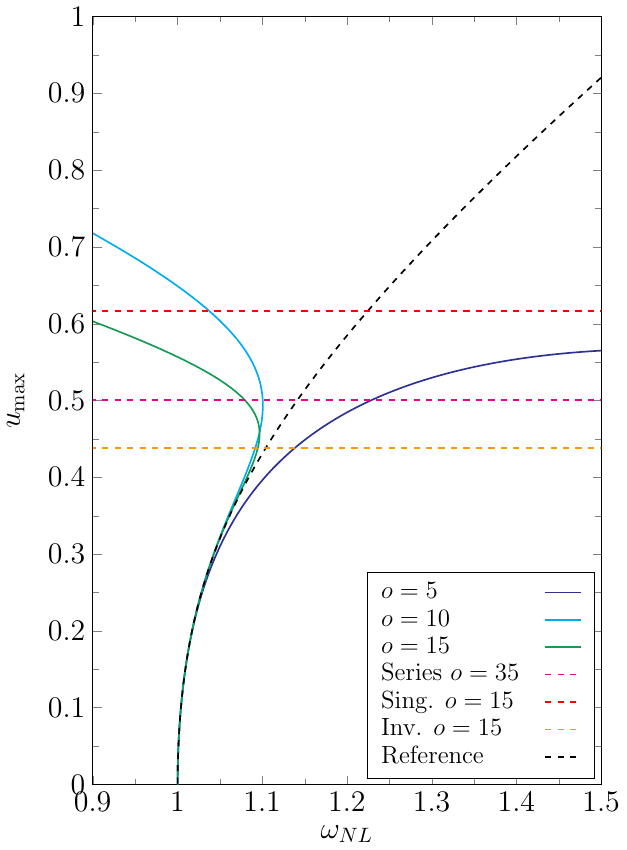}
		\caption{$\omega_2 = 0.637$}
		\label{fig:2 DOF cubic conservative unforced - omega2 = 0.637 - Backbones}
	\end{subfigure}
	\caption{Backbone curves for the two-DOF system. The reference solution (dashed black line) obtained by numerical continuation (using the package matcont \cite{dhooge2004matcont}) is compared to ROMs using with a CNF style with increasing orders $o$. The different validity limits criteria are reported with horizontal dashed colour lines. Parameter values are set as $\omega_1 = h^1_{111} = h^2_{111} = 1.0$.}
	\label{fig:2 DOF cubic conservative unforced - Backbones}
\end{figure}

\begin{figure}[h]
	\centering
	\begin{subfigure}{0.49\textwidth}
		\includegraphics[width=\textwidth]{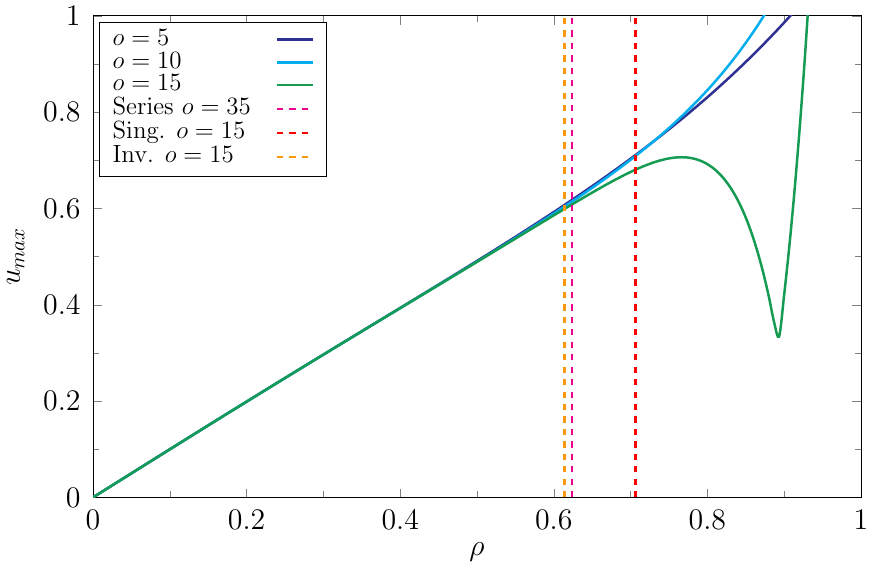}
		\caption{$\omega_2 = 1.57$}
		\label{fig:2 DOF cubic conservative unforced - omega2 = 1.57 - Max disp}
	\end{subfigure}
	\begin{subfigure}{0.49\textwidth}
		\includegraphics[width=\textwidth]{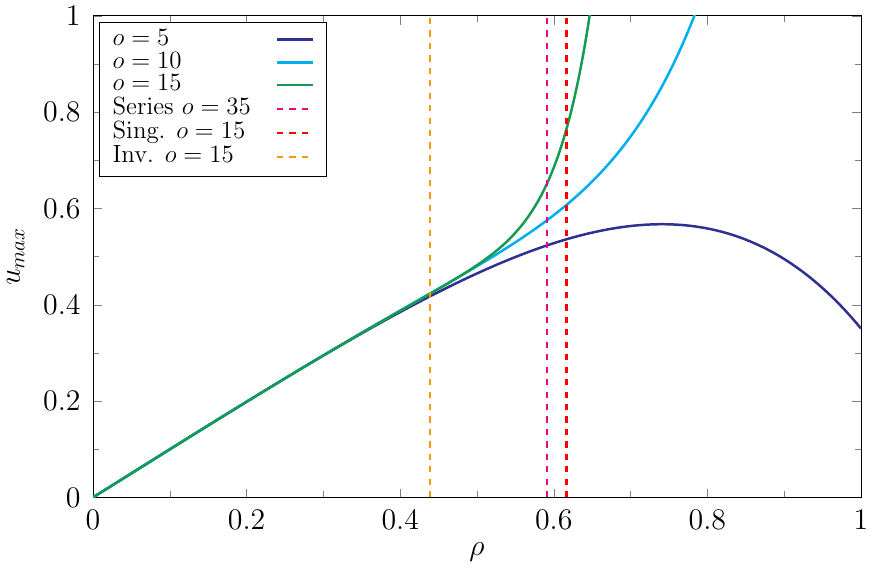}
		\caption{$\omega_2 = 0.637$}
		\label{fig:2 DOF cubic conservative unforced - omega2 = 0.637 - Max disp}
	\end{subfigure}
	\caption{Relation between amplitude of the normal variables $\rho$ and maximum physical displacement amplitude $u_{\text{max}}$ for the unforced and undamped 2-DOF oscillator. The validity limits obtained by the proposed criteria are also presented, in dashed lines. Parameter values are fixed as $\omega_1 = h^1_{111} = h^2_{111} = 1$.}
	\label{fig:2 DOF cubic conservative unforced - Max disp}
\end{figure}

Summaries with the results from all criteria are given in \cref{tab:2 DOF cubic conservative unforced - omega2 = 0.637 - Rho values,tab:2 DOF cubic conservative unforced - omega2 = 1.57 - Rho values}, and their values are graphically compared with backbone curves and with a reference obtained by numerical continuation in \cref{fig:2 DOF cubic conservative unforced - Backbones}. The same trend as observed before remains: the invariance equation criterion with $\varepsilon=0.01$ serves as a lower bound for the zone where graphically the parametrisation seems to stop converging, while the series and singularity criteria give an upper estimate. For \cref{fig:2 DOF cubic conservative unforced - omega2 = 0.637 - Backbones}, though, the region indicated by the criteria is quite extensive, and includes zones where the backbone has completely diverged from the reference solution. This behaviour can be understood by looking at the relationship between normal variables amplitude $\rho$ and maximum displacement amplitude $u_{\text{max}}$, given in \cref{fig:2 DOF cubic conservative unforced - Max disp}. It can be noted that for the case $\omega_2 = 0.637$, the validity limits obtained by the singularity and invariance equation criteria enter a zone where the relationship $\rho - u_{\text{max}}$ becomes excessively steep, whereas this does not happen for the case $\omega_2 = 1.57$. The commented behaviour indicates a certain difficulty of the singularity and series criteria to deal with the excessive growth in the $\rho-u_{\text{max}}$ for high degrees of expansions. A large gap between the validity limits indicated by the different criteria should also be a warning indicator for the analyst of such a phenomenon. In order to mitigate this effect, a lower order expansion can be chosen only to compute the maximum physical displacement for a given value of reduced variables amplitude. For instance, if $o=10$ is selected, \cref{fig:2 DOF cubic conservative unforced - omega2 = 0.637 - Max disp} shows that the region where convergence is lost is much smaller, aligning with qualitative expectations. However, it should be noted that in all conducted analyses, this phenomenon was observed only at high expansion orders, which are rarely used in practice.

\begin{table}[h]
	\npdecimalsign{.}
	\nprounddigits{3}
	\centering
	\begin{tabular}{c|c|c|c|c}
		Criterion & Cauchy & \text{d'Alembert} & Singularity & Invariance \\ \hline
		$\rho$ & 0.730 & 0.639 & 0.766 & 0.630 \\
		$u_{\text{max}}$ & 0.696 & 0.623 & 0.706 & 0.614 \\
	\end{tabular}
	\caption{Validity limits for the unforced and undamped 2-DOF oscillator with $\omega_2 = 1.57$. Values are given in both the amplitude of the normal coordinates and maximum physical displacement of the oscillator. Orders of expansion are 15 for the singularity and invariance equation criteria and 35 for the series criteria. The tolerance for the invariance equation is chosen as $\varepsilon = 1\%$. Parameter values are fixed as $\omega_1 = h^1_{111} = h^2_{111} = 1.0$.}
	\label{tab:2 DOF cubic conservative unforced - omega2 = 1.57 - Rho values}
\end{table}

\begin{table}[h]
	\npdecimalsign{.}
	\nprounddigits{3}
	\centering
	\begin{tabular}{c|c|c|c|c}
		Criterion & Cauchy & \text{d'Alembert} & Singularity & Invariance \\ \hline
		$\rho$ & 0.517 & - & 0.580 & 0.456 \\
		$u_{\text{max}}$ & 0.501 & - & 0.617 & 0.439 \\
	\end{tabular}
	\caption{Validity limits for the unforced and undamped 2-DOF oscillator with $\omega_2 = 0.637$. Values are given in both the amplitude of the normal coordinates and maximum physical displacement of the oscillator. Orders of expansion are 15 for the singularity and invariance equation criteria and 35 for the series criteria. The tolerance for the invariance equation is chosen as $\varepsilon = 1\%$. Parameter values are fixed as $\omega_1 = h^1_{111} = h^2_{111} = 1.0$.}
	\label{tab:2 DOF cubic conservative unforced - omega2 = 0.637 - Rho values}
\end{table}

The forced-damped case is now investigated. The focus is set on the case where $\omega_2 = 1.57$. The damping coefficients are selected as: $\xi_1=\xi_2=0.05$. Primary resonance is studied, and the expansion forcing frequency is chosen as $\Omega = \omega_1$. For this case, the validity limit obtained by the invariance equation criterion for the unforced and undamped case is $\rho = 0.63$. In order to determine the adequate level of forcing for this example, we employ the procedure described in \cref{sec:a_priori_amplitude}. In this case, by inspecting the reduced dynamics coefficients, we remark that $f_3^R = 0$ and thus we can not employ \cref{eq:kappa_max_gen}. Instead, we apply \cref{eq:max_kappa_fI}, and find that the maximum forcing amplitude that the system can be subjected to, associated to the radius of convergence obtained by the invariance equation criterion, is $\kappa = 0.0629$. Since the invariance equation usually gives a lower bound, we will instead choose a value of $\kappa = 0.07$, approximately $10\%$ higher than the computed one.

\begin{figure}[h]
	\centering
	\includegraphics[width=0.5\textwidth]{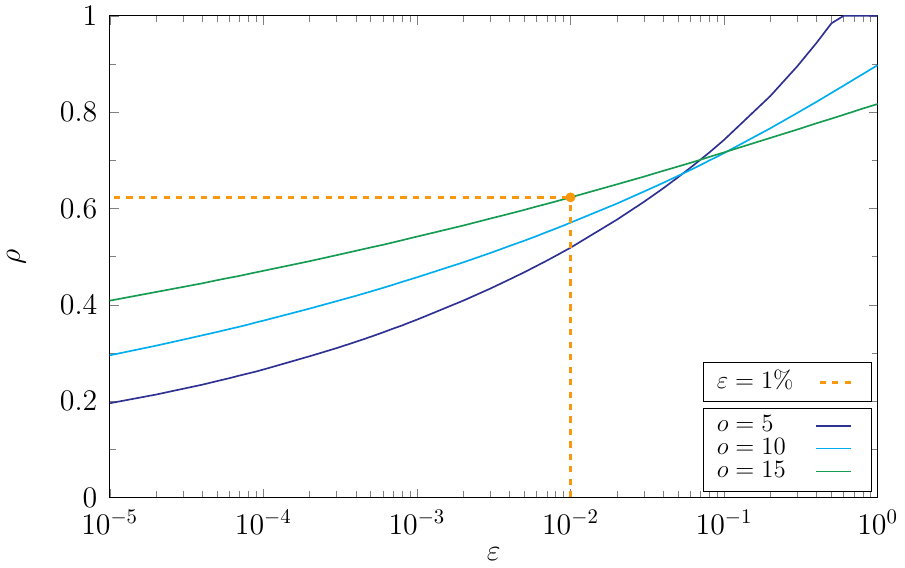}
	\caption{Evolution of the radius of convergence with the tolerance value for the 2-DOF oscillator in primary resonance using the invariance equation criterion. Parameter values are fixed as $\omega_1 = h^1_{111} = h^2_{111} = 1$, $\xi_1=\xi_2=0.05$ and $\kappa = 0.07$.}
	\label{fig:2 DOF cubic damped primary - omega2 = 1.57 - Invariance equation criterion}
\end{figure}

For this level of forcing, $\kappa=0.07$ and we present the curves associated with the invariance equation criterion in \cref{fig:2 DOF cubic damped primary - omega2 = 1.57 - Invariance equation criterion}. Visually, it can be seen that for a tolerance value of $\varepsilon = 1\%$, the value of the convergence radius does seem to be around $\rho = 0.63$, the value for the unforced system, as we have assumed in order to estimate $\kappa$.

\begin{figure}[h]
\begin{minipage}[t]{.48\textwidth}
	\centering
	\includegraphics[width=\textwidth,valign=t]{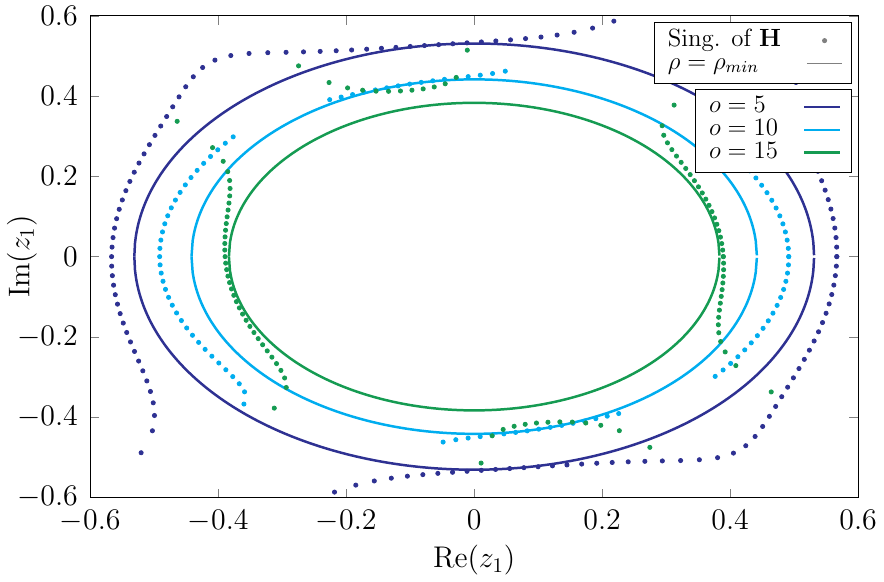}
	\captionof{figure}{Validity limit calculated by the criterion of the singularity of the homological operator for the 2-DOF oscillator in primary resonance as a function of the parametrisation order. The lines represent circles of radius equal to the minimum amplitude among points of the same order. Parameter values are fixed as $\omega_1 = h^1_{111} = h^2_{111} = 1.0$, $\xi_1=\xi_2=0.05$ and $\kappa = 0.07$.}
	\label{fig:2 DOF cubic damped primary - omega2 = 1.57 - Singularity criterion}
\end{minipage}%
\hspace{0.04\textwidth}
\begin{minipage}[t]{.48\textwidth}
	\centering
	\includegraphics[width=\textwidth,valign=t]{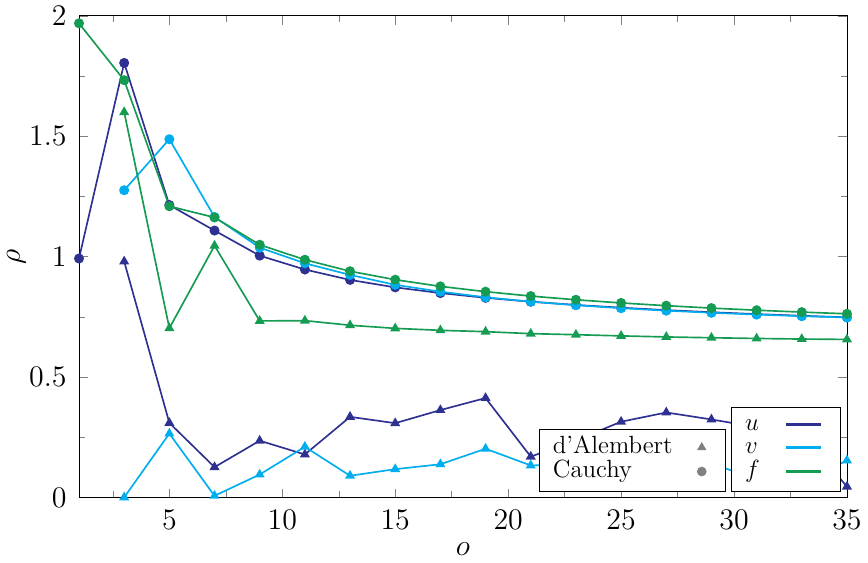}
	\captionof{figure}{Validity limit calculated by the series criteria for the 2-DOF oscillator in primary resonance as a function of the parametrisation order. Parameter values are fixed as $\omega_1 = h^1_{111} = h^2_{111} = 1.0$, $\xi_1=\xi_2=0.05$ and $\kappa = 0.07$.}
	\label{fig:2 DOF cubic damped primary - omega2 = 1.57 - Series criteria}
\end{minipage}
\end{figure}

\begin{figure}[h]
	\centering
	\includegraphics[width=0.5\textwidth]{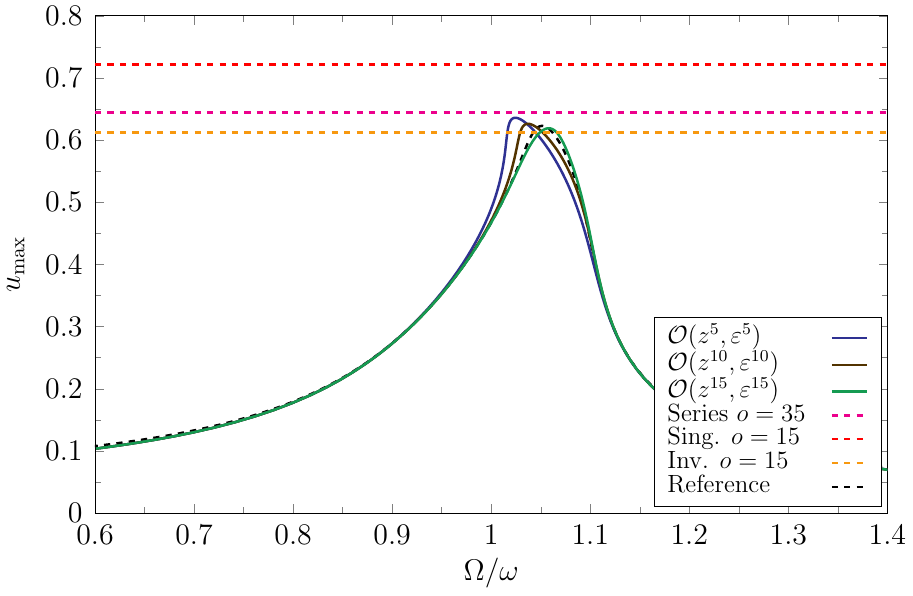}
	\caption{Frequency response curves for the 2-DOF oscillator in primary resonance. Different orders of parametrisation are compared with a reference solution computed by numerical continuation employing package matcont \cite{dhooge2004matcont}. The validity limits calculated by the different proposed convergence criteria are also shown. Parameter values are fixed as $\omega_1 = h^1_{111} = h^2_{111} = 1.0$, $\xi_1=\xi_2=0.05$ and $\kappa = 0.07$.}
	\label{fig:2 DOF cubic damped primary - omega2 = 1.57 - FRCs}
\end{figure}

We continue by examining the singularity criterion. The results obtained for this case are given in \cref{fig:2 DOF cubic damped primary - omega2 = 1.57 - Singularity criterion}. Qualitatively, they are close to the unforced case, but the effect of the excitation is revealed in the asymmetric distribution of singular points. Finally, the series criteria are also presented, in \cref{fig:2 DOF cubic damped primary - omega2 = 1.57 - Series criteria}. As usual, Cauchy's rule presents a better convergence than d'Alembert's. However, the lowest radius of convergence value is the one given by the former criterion. The various radii of convergence are summarised in \cref{tab:2 DOF cubic damped primary - omega2 = 1.57 - Rho values}, and the FRCs obtained for this case are presented in \cref{fig:2 DOF cubic damped primary - omega2 = 1.57 - FRCs}. It can be seen that the peak of the FRCs coincides almost exactly with the predicted validity limit, which validates the proposed procedure to estimate $\kappa$. Also, although an order 15 parametrisation was computed, the ROM and FOM curves start to slightly disagree at the peak of the FRCs, indicating that the validity limit is indeed reached. All of the validity limits present a good predictive behaviour in this case, and the zone indicated by the invariance and series criteria for loss of convergence is quite narrow.

\begin{table}[h]
	\npdecimalsign{.}
	\nprounddigits{3}
	\centering
	\begin{tabular}{c|c|c|c|c}
		Criterion & Cauchy & \text{d'Alembert} & Singularity & Invariance \\ \hline
		$\rho$ & 0.748 & 0.656 & 0.766 & 0.623 \\
		$u_{\text{max}}$ & 0.720 & 0.645 & 0.712 & 0.613 \\
	\end{tabular}
	\caption{Validity limits for the 2-DOF oscillator in primary resonance. Values are given in both the amplitude of the normal coordinates and maximum physical displacement of the oscillator. Orders of expansion are 15 for the singularity and invariance equation criteria and 35 for the series criteria. The tolerance for the invariance equation is chosen as $\varepsilon = 1\%$. Parameter values are fixed as $\omega_1 = h^1_{111} = h^2_{111} = 1.0$, $\xi_1=\xi_2=0.05$ and $\kappa = 0.07$}
	\label{tab:2 DOF cubic damped primary - omega2 = 1.57 - Rho values}
\end{table}

The cases studied in this section underline that an effective procedure can be set to derive the validity limits of the forced case from the unforced analysis. The next section considers a large system stemming from the finite element discretisation of a continuous structure, to assess the different criteria and the proposed rules in a more adverse situation.

\subsection{Continuous systems: a cantilever beam discretised with finite elements} \label{sec:Beam}
In this section, we consider a cantilever beam modelled with 2D quadratic finite elements and a plane stress assumption. The beam has a length of $L=1\text{ m}$ and an in-plane thickness of $H = 0.02\text{ m}$. Although large displacements are considered, strains are supposed to be small due to the small thickness of the beam, and a material following Saint-Venant Kirchhoff's law, having Young's modulus of $E = 104 \cdot 10^9\text{ GPa}$, Poisson's coefficient of $\nu = 0.3$ and density of $\rho = 4400\text{ kg}/\text{m}^3$, is chosen. Such a cantilever beam has already been studied in~\cite{vizza21high} with three-dimensional elements, as well as in~\cite{Grolet2025} with a beam model. Besides, the 2D model studied here has already been used in~\cite{Stabile:follow}, and the interested reader is referred to that article for more details on the finite element formulation.

The unforced and undamped response of the system is first studied. The invariant manifold corresponding to the first bending mode for this example presents a particular feature: a folding point is present for a moderate value of the vibration amplitude~\cite{vizza21high}. Consequently, only normal form styles of parametrisation are adequate to treat this example. Using a graph style leads to a rapid failure of the computed NNM in that case, because of the impossibility of passing over the fold. This aspect should be reflected by the calculated validity limits, and thus, differently from the previous sections, both the CNF and graph styles are investigated.

\begin{figure}[h]
	\centering
	\begin{subfigure}{0.49\textwidth}
		\includegraphics[width=\textwidth]{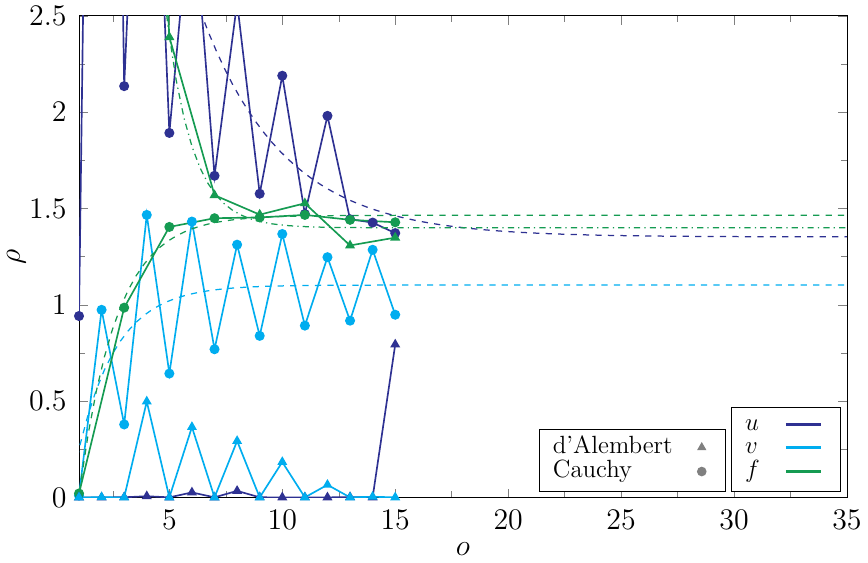}
		\caption{CNF style}
		\label{fig:Cantilever unforced - CNF style - Series criteria}
	\end{subfigure}
	\begin{subfigure}{0.49\textwidth}
		\includegraphics[width=\textwidth]{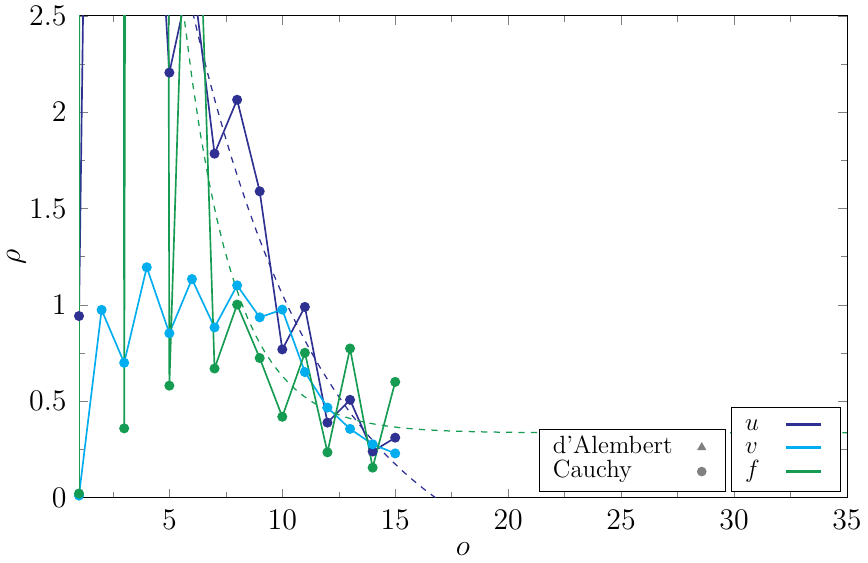}
		\caption{Graph style}
	\end{subfigure}
	\caption{Validity limit calculated by the series criteria for the unforced cantilever beam as a function of the parametrisation order. For each of the criteria, the displacement, velocity and reduced dynamics power series are studied. Results obtained by extrapolation of parametrisations up to order 15 are also presented, in dashed lines.}
	\label{fig:Cantilever unforced - Series criteria}
\end{figure}

The results obtained by the series criteria are shown in \cref{fig:Cantilever unforced - Series criteria}. In this example, and for structures discretised by finite elements in general, it is generally too costly to expand the parametrisation to such high orders as in the previous cases, and thus the convergence radius of the method must be extrapolated from only a few points corresponding to low-order computations. Specifically, the dashed lines in the figure correspond to curves obtained by doing a least squares fit of the values up to order $o = 15$, considering functions of exponential decay. Such fits are displayed only for curves for which a convergent trend was observed, mainly the Cauchy criterion ones. Besides, due to the oscillatory behaviour from one order to another, the fit is performed by considering the mean values of successive data points.

The curves corresponding to a CNF style visually present an acceptable convergence, with the exponential extrapolation yielding sensible values for the convergence radius. For the graph style, however, convergence is much worse. As a matter of fact, for the d'Alembert criterion, the estimated convergence radii oscillate between extremely high values and values close to zero, such that the associated curves were omitted from the figure. The best result for the extrapolation procedure was associated with the reduced dynamics curve for the Cauchy criterion. Nevertheless, even though the convergence radius in $\rho$ is small, it corresponds to large values in $u_{\text{max}}$, indicating that the convergence limit has already been reached. Thus, the retained value of $\rho$ for the series criteria with the graph style parametrisation was the lowest computed at order 15, corresponding to the Cauchy criterion on velocity. Note that, in that specific case, the extrapolation did not provide meaningful values, and the last computed point has thus been used without resorting to an extrapolation.

\begin{figure}[h]
\begin{minipage}[t]{.48\textwidth}
	\centering
	\includegraphics[width=\textwidth,valign=t]{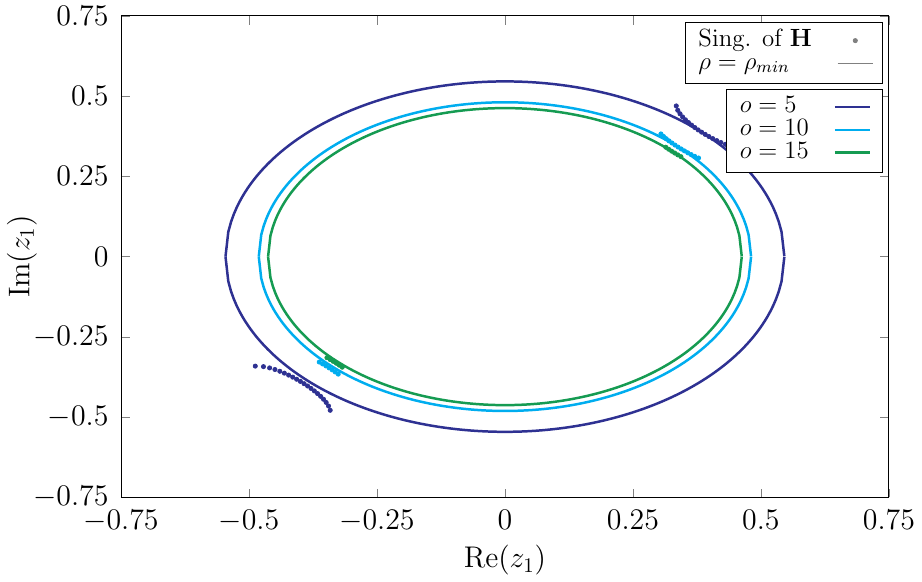}
	\captionof{figure}{Validity limit calculated by the criterion of the singularity of the homological operator for the unforced cantilever beam as a function of the parametrisation order. The lines represent circles of radius equal to the minimum amplitude among points of the same order.}
	\label{fig:Cantilever unforced - CNF style - Singularity criterion}
\end{minipage}
\hspace{0.04\textwidth}
\begin{minipage}[t]{.48\textwidth}
	\centering
	\includegraphics[width=\textwidth,valign=t]{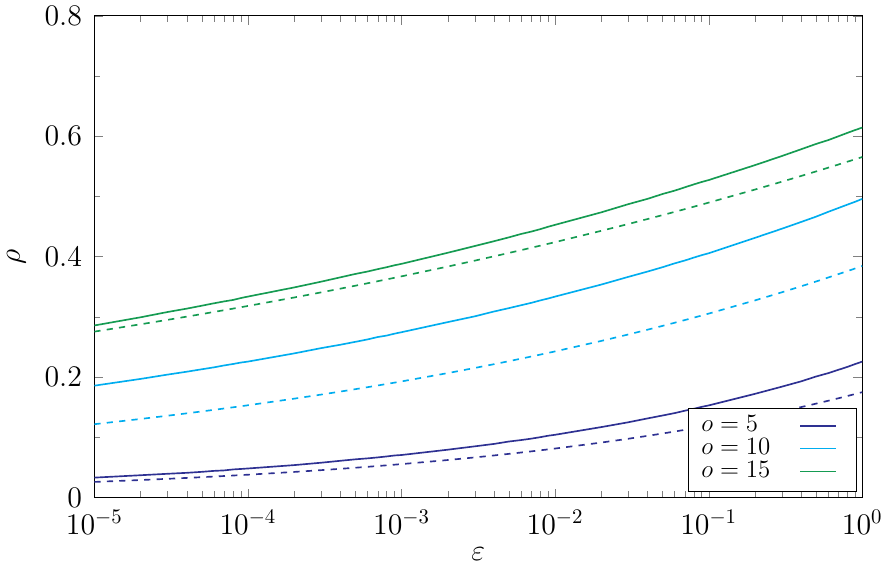}
	\captionof{figure}{Evolution of the radius of convergence with the tolerance value for the unforced cantilever beam using the invariance equation criterion. The continuous lines indicate results obtained by \cref{eq:inv_eq_error}, while the dashed ones represent the simplification introduced in \cref{eq:inv_eq_error_ANM}.}
	\label{fig:Cantilever unforced - CNF style - Invariance equation criterion}
\end{minipage}
\end{figure}

For the singularity criterion, results are presented in~\cref{fig:Cantilever unforced - CNF style - Singularity criterion} for the CNF style. A key observation is that, unlike in the other examples, this situation yielded only a limited number of singularity points. At first, one might think that this could be an effect of the choice of the submatrix used in the calculation of the determinant of the homological operator. Specifically, vertical displacement and velocity degrees of freedom of the point on the lower side and on the tip of the beam were selected. However, a more general approach based on singular values for this example has been tested in \cref{sec:sing_comp}, and although this approach finds more points where the operator's rank changes, the results in terms of validity limits are the same. This indicates that there are indeed no singularities, and that the submatrix procedure can be safely used. For graph style, it was not possible to find singularities for the homological operator. This seems sensible in some sense, as the relationship between normal coordinates and physical variables is this case is a functional one. Concerning the invariance equation criterion, the evolution of the convergence radius with the tolerance $\varepsilon$ is depicted in \cref{fig:Cantilever unforced - CNF style - Invariance equation criterion} for the CNF style of parametrisation, with the simplification brought about by \cref{eq:inv_eq_error_ANM} in dashed line. The simplified computation yields excellent results for odd orders of parametrisations, but higher discrepancies can be verified for even orders.

\begin{figure}[h]
	\centering
	\begin{subfigure}{0.49\textwidth}
		\includegraphics[width=\textwidth]{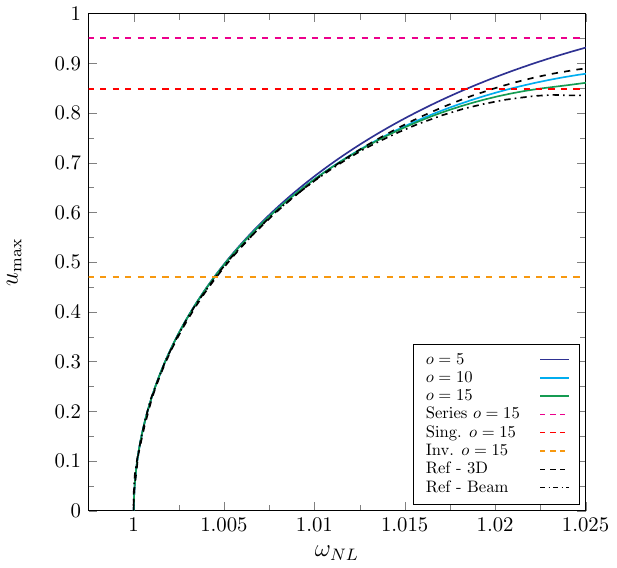}
		\caption{CNF style}
		\label{fig:Cantilever unforced - CNF style - Backbones}
	\end{subfigure}
	\begin{subfigure}{0.49\textwidth}
		\includegraphics[width=\textwidth]{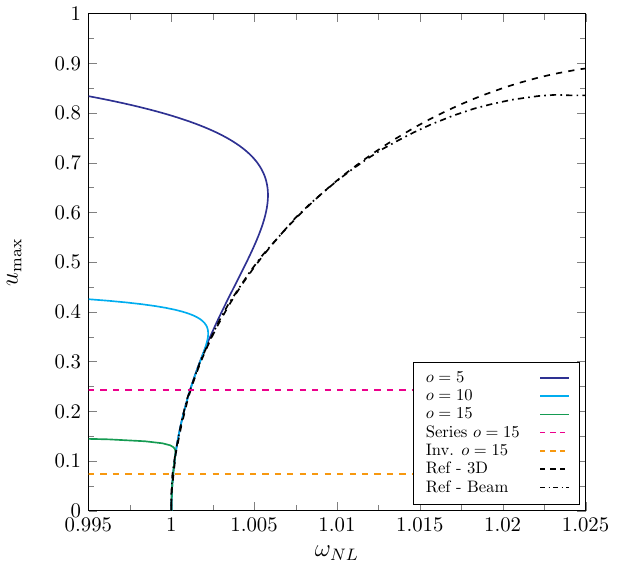}
		\caption{Graph style}
		\label{fig:Cantilever unforced - Graph style - Backbones}
	\end{subfigure}
	\caption{Comparison of the backbone curves given by different orders of parametrisation with two reference solutions: one obtained from a high order parametrisation applied to a 3D FE model \cite{vizza21high} and another by numerical continuation of a geometrically exact beam model for the cantilever beam \cite{Grolet2025}. The different proposed validity criteria are shown in dashed colored lines.}
	\label{fig:Cantilever unforced - Backbones}
\end{figure}

To conclude, the obtained convergence radii are summarised in \cref{tab:Cantilever unforced - CNF style - Rho values,tab:Cantilever unforced - Graph style - Rho values} and visually depicted, together with the backbone curves with increasing orders, in \cref{fig:Cantilever unforced - Backbones}. In addition to the calculated backbones, reference solutions taken from \cite{vizza21high,Grolet2025} are also shown. The solution of \cite{vizza21high} was obtained by a high-order parametrisation method with a 3D finite element formulation, while that presented in \cite{Grolet2025} stems from a harmonic balance calculation on a geometrically exact beam model implementing Timoshenko kinematics. Considering first the CNF style, the zone between the invariance equation and the other two criteria is quite extensive in this case. Yet, it is important to remember that the backbone only summarises information about a single displacement, and thus the relatively low value of the limit obtained by the invariance equation in this situation can indicate that in another part of the manifold the parametrisation starts to fail. Additionally, it is interesting to point out that the way the ROM solutions fail in this example is different from all of the other test cases, with a smooth departure from the full-order solution. This could be a cause for the relatively large zone of divergence predicted by the validity limits. Now, looking at the results for the graph style, the parametrisation almost immediately fails, due to a fold in the manifold of the cantilever beam. This is clearly reflected in the validity limits, which are correspondingly low and accurately indicate the region where the parametrisation breaks down.

\begin{table}[h]
	\npdecimalsign{.}
	\nprounddigits{3}
	\centering
	\begin{tabular}{c|c|c|c|c}
		Criterion & Cauchy & \text{d'Alembert} & Singularity & Invariance \\ \hline
		$\rho$ & 1.102 & 1.400 & 0.925 & 0.453 \\
		$u_{\text{max}}$ & 0.951 & 2.476 & 0.849 & 0.473 \\
	\end{tabular}
	\caption{Validity limits for the unforced cantilever beam studied with complex normal form style. Values are given in both the amplitude of the normal coordinates and maximum physical displacement of the oscillator. Orders of expansion are 15 for the singularity and invariance equation criteria, and the values for the series criteria are obtained by extrapolation. The tolerance for the invariance equation is chosen as $\varepsilon = 1\%$.}
	\label{tab:Cantilever unforced - CNF style - Rho values}
\end{table}

\begin{table}[h]
	\npdecimalsign{.}
	\nprounddigits{3}
	\centering
	\begin{tabular}{c|c|c|c|c}
		Criterion & Cauchy & \text{d'Alembert} & Singularity & Invariance \\ \hline
		$\rho$ & 0.228 & - & - & 0.07 \\
		$u_{\text{max}}$ & 0.243 & - & - & 0.0744 \\
	\end{tabular}
	\caption{Validity limits for the unforced cantilever beam studied with graph style. Values are given in both the amplitude of the normal coordinates and maximum physical displacement of the oscillator. Orders of expansion are 15 for the singularity and invariance equation criteria, and the values for the series criteria are obtained by extrapolation. The tolerance for the invariance equation is chosen as $\varepsilon = 1\%$.}
	\label{tab:Cantilever unforced - Graph style - Rho values}
\end{table}

After studying the system in free vibration, we consider now a case with forcing and damping. Specifically, we choose a forcing that is proportional to the first bending mode, such that
\begin{equation}
	\bfF_c = \kappa \bfM \bfphi_{B1},
\end{equation}
where $\bfM$ denotes the mass matrix of the system and $\bfphi_{B1}$ the modal shape (consisting of only displacements) of the first bending mode, with a mass normalisation condition. The expansion frequency for the parametrisation is chosen as the first bending frequency: $\Omega=98.714\text{ rad/s}$, and the damping is mass proportional, with a proportionality coefficient of $98.714/1000$, corresponding to a damping ratio of $\xi=0.0005$. Only the CNF style is considered. For this case, the validity limit obtained by the invariance equation criterion for the unforced and undamped case is $\rho = 0.925$. The applied forcing level is chosen once again by means of the procedure described in \cref{sec:a_priori_amplitude}. Employing \cref{eq:kappa_max_gen}, we find the corresponding forcing amplitude to be approximately $\kappa = 7 \text{ m/s$^2$}$. For this example, we only summarise the validity limits obtained by the three different criteria in \cref{tab:Cantilever damped primary - Rho values}. Specifically, we have computed parametrisations only up to order 10, as, due to the inclusion of the forcing and the number of degrees of freedom, going further implies important computation times.

\begin{table}[h]
	\npdecimalsign{.}
	\nprounddigits{3}
	\centering
	\begin{tabular}{c|c|c|c|c}
		Criterion & Cauchy & \text{d'Alembert} & Singularity & Invariance \\ \hline
		$\rho$ & 0.948 & 1.420 & 0.961 & 0.276 \\
		$u_{\text{max}}$ & 0.871 & 1.414 & 0.876 & 0.292 \\
	\end{tabular}
	\caption{Validity limits for the cantilever beam in primary resonance. Values are given in both the amplitude of the normal coordinates and maximum physical displacement of the oscillator. The criteria were calculated at an expansion order equal to 10, with the exponential extrapolation procedure used for the series criteria. The tolerance for the invariance equation is chosen as $\varepsilon = 1\%$.}
	\label{tab:Cantilever damped primary - Rho values}
\end{table}

\begin{figure}[h]
	\centering
	\includegraphics[width=0.5\textwidth]{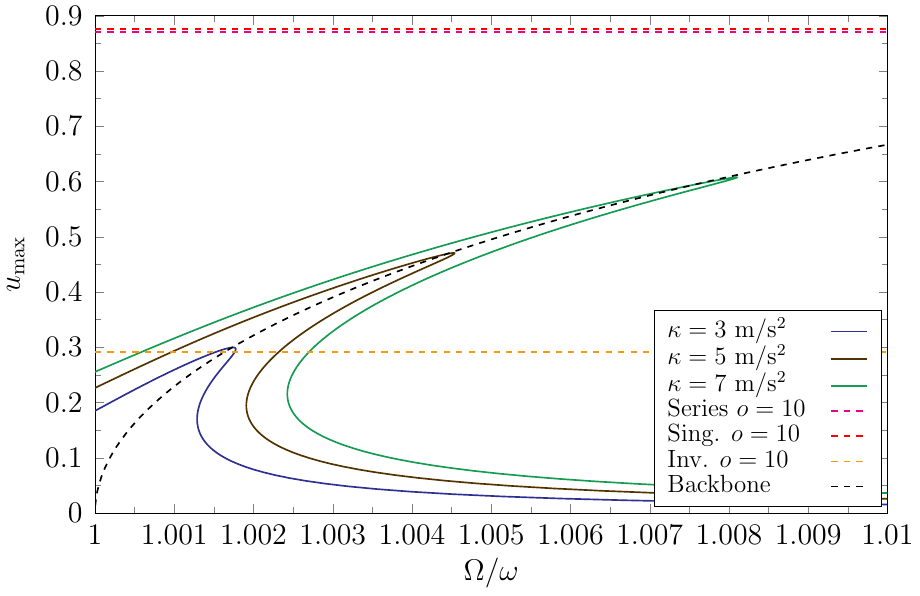}
	\caption{Frequency response curves for the cantilever beam in primary resonance. Different orders of parametrisation and the validity limits calculated by the different proposed convergence criteria are shown.}
	\label{fig:Cantilever damped forced - FRCs}
\end{figure}

The frequency response curves for this case with parametrisation order equal to 10 and $\kappa = 3,5,7 \text{ m/s$^2$}$ are given in \cref{fig:Cantilever damped forced - FRCs}, together with a graphical depiction of the validity limits and the backbone curve of the problem. In this example we do not compute a full order solution, as would be the case in a real-life application of the proposed methodology. It is known, however, from a previous study \cite{opreni22high}, that the ROM converges up to $\kappa = 5 \text{ m/s$^2$}$, and thus the predicted value of $\kappa = 7 \text{ m/s$^2$}$ is in line with the expected order of magnitude for this case. From the figure, it can be seen that, as was the case for the unforced scenario, the validity criteria indicate a large zone of loss of convergence for the reduced order model. Nevertheless, it is a region that does encompass the peak of the FRC for $\kappa = 7 \text{ m/s$^2$}$.

\section{Practical assessment of the criteria and guidelines for the analysis}

After illustrating the application of the different validity limit criteria through several examples, this section turns to a practical comparison of their performance and provides guidance on how to employ them in routine analyses. We begin by summarizing the main characteristics observed for each criterion across the examples.

The invariance equation criterion appears to be the most effective for estimating the expected convergence range of the approximation. It was shown to be robust, requiring no adjustments for its application to large models, while consistently yielding a coherent value for the radius of convergence. However, it was also the most time-consuming criterion to compute, and its main drawback lies in the fact that its results depend on the choice of a tolerance value, which is not always straightforward. For the presented examples, the value $\varepsilon = 10^{-2}$ yielded generally a good approximation of the amplitude for which the approximated results depart from the reference. When a single mode is considered and the system is autonomous, it is possible to apply the simplification given by \cref{eq:inv_eq_error_ANM}, which gives precise results when compared to the criterion in its full form and drastically reduces its computation time.  Across all the examples presented, this criterion established a lower bound on the validity limit, though this result is contingent upon the tolerance chosen for the invariance equation error.

The singularity criterion provides an upper bound for the convergence radius of the approximation and was generally accurate. However, for the case of the forced beam, only a few singularity points were found. As both the strategy of computing the determinant of a submatrix of the homological operator and tracking its singular values were tested, we conclude that this is a feature of the selected example. A more robust approach to perform this calculation may exist, but this remains an open question. The criterion should be tested in other high-dimensional systems to better evaluate its performance in adverse conditions.

Considering now the series criteria, an upper bound for the validity limit is also found. For almost all of the examples, Cauchy's criterion worked better than d'Alembert's, and as a general rule, we would advise the use of the former. Their computation is in general fast, provided a high enough parametrisation has already been computed. Since they presuppose convergence of the power series stemming from the method, this can be limiting, especially for high-dimensional systems, but the extrapolation procedure proposed seems to yield good results and to make their use possible.

Finally, for practical analyses, we suggest that even when a forced system is considered, the underlying unforced and conservative system should also be examined. This enables the computation of the backbone curve and allows one to estimate the maximum admissible forcing value using the procedure described in~\cref{sec:a_priori_amplitude}.
Concerning the actual computation of the validity limit, we suggest computing the three criteria: the invariance equation provides a lower bound, while the minimum value between the singularity and series criteria can be taken as an upper bound for the convergence radius. 
This procedure identifies three distinct amplitude zones. When the vibration amplitude falls below the lower bound of the validity limit, the ROM may be safely applied; when it exceeds the upper bound, the ROM is almost certainly invalid. In the intermediate region, however, the decision must be made by the analyst, with appropriate validation through complementary approaches.

\section{Conclusion} \label{sec:Conclusions}

The direct parametrisation method for invariant manifold is an efficient reduction method for vibrating systems, which relies on a local theory. Thanks to automated arbitrary order expansions, the derived ROMs provide exact results as long as the validity limit of the expansions are not met. In this context, the aim of this contribution is to assess different criteria proposed in the literature, for the specific case of nonlinear vibrating systems. Reduction to a single NNM is also assumed for this investigation, and practical guidelines have been given for determining the convergence zone of the issuing ROMs without the need to compute full-order solutions. Different cases have been showcased with increasing complexity: conservative systems (backbone curves), forced response in primary and superharmonic resonance, from single DOF problem to finite element structures.

The main outcomes can be summarised as follows. Throughout the examples, it has been found that the error in the invariance equation with a tolerance $\varepsilon=10^{-2}$, provides an accurate lower bound. An alternative, simplified version of this criterion has been proposed, which provides a much faster way to compute the validity limit, but only for single DOF problems without forcing. It has been found that this simplified criteria on the invariance equation also gives an accurate and easy-to-use lower bound. An adaptation of the criterion proposed in~\cite{LamarqueUP} has been proposed for application to the DPIM framework. This singularity criterion, together with the series criteria proposed in~\cite{Grolet2025}, has been shown to provide an upper bound in practical examples. For the series criteria, numerical results demonstrated that extrapolations can be used quite safely from a small number of low-order computations, to infer the limit value. It has been thus proposed, in practice, to use the three criteria, and use the interval in between the provided bounds as a critical region where care about the approximation should be paid. Finally, it has been shown that the validity limit for the forced case can be derived from an analysis of the unforced problem, by using the maximum amplitude found and a simple rule that predicts the maximum forcing amplitude allowed. Examples showed that this provides a robust and accurate workout for practical analysis.

To conclude, it is mentioned that more tests on a larger number of continuous systems should give a better view of the behaviour of the criteria for predicting the validity limits. It is thus advocated to use them in the numerical codes performing the parametrisation method for vibrating systems. Also, extension to the cases where more than one master mode is selected should be further investigated to see if simple and accurate criteria could also be used. 

\section*{Acknowledgments}
The authors would like to thank Olivier Thomas and Claude-Henri Lamarque, with whom numerous in-depth discussions were held on issues related to the limits of validity, whether in normal form theory or in the parametrisation method for invariant manifolds. They are gratefully acknowledged for the many discussions and ideas that have greatly contributed to our reflections. We also thank Maxime Breden, for the fruitful discussion on the effect of the normalisation of the eigenvectors on the validity limits of the parametrisation method. His aid and the time he has spent with us are greatly appreciated.

\section*{Funding}
The work received no additional funding.

\section*{Conflict of interest} 
The authors declare that they have no conflict of interest.

\section*{Data availability statement}
The codes developed and data used for the presented analyses are available upon reasonable request to the authors.

\section*{Author contributions}
Conceptualisation: all authors contributed equally; Methodology: all authors contributed equally; Formal analysis and investigation: Andr\'e de F. Stabile; Software development: Andr\'e de F. Stabile, Aur\'elien Grolet; Writing - original draft preparation: Andr\'e de F. Stabile, Cyril Touz\'e; Writing - review and editing: Andr\'e de F. Stabile, Cyril Touz\'e, Aur\'elien Grolet, Alessandra Vizzaccaro;
Supervision: Cyril Touz\'e, Alessandra Vizzaccaro, Aur\'elien Grolet; 

\bibliographystyle{spmpsci}      % mathematics and physical sciences
\bibliography{biblioROM}

\appendix

\counterwithin{figure}{section}
\renewcommand\thefigure{\thesection\arabic{figure}}

\section{Simplification of the invariance equation criterion} \label{sec:inv_eq_simpl}

In this appendix, we give further details on the approximation of the numerator of~\cref{eq:inv_eq_error}, which is related to the error in the invariance equation $\bfE$ defined in~\cref{eq:numinvaerror}. This quantity can be generically decomposed into its monomial components by
\begin{equation}
	\bfE = \sum_{p=1}^{\infty} \sum_{k=1}^{m_p} \bfE^{(p,k)} \bfz^{\alphavec(p,k)},
\end{equation}
where no truncation of the power series has been performed. In \cite{Vizza:superDPIM}, it is shown that, for a fixed monomial, the term $\bfE^{(p,k)}$ can be written as
\begin{equation}
	\bfE^{(p,k)} = \left( \sigma^{(p,k)} \bfB - \bfA \right) \bfW^{(p,k)} + \sum_{s=1}^{d} \bfB \bfY_s f_s^{(p,k)} - \bfR^{(p,k)}.
\end{equation}
This expression also corresponds to~\cref{eq:homological}, which is rewritten in an augmented formulation to deal with the singular cases corresponding to resonances.

Let us now consider a parametrisation that has been computed up to order $o$. This means that $\forall p \leq o, \bfE^{(p,k)} = \bf0$, since $\bfW^{(p,k)}$ and $\bff^{(p,k)}$ are obtained as solutions of \cref{eq:homological}. Furthermore, in this situation, all of the nonlinear mapping and reduced dynamics coefficients for orders higher than $o$ vanish. In particular, at order $o+1$, one simply obtains:
\begin{equation}
	\bfE^{(o+1,k)} = \bfR^{(o+1,k)}.
\end{equation}
Therefore, by approximating $\bfE$ by its leading-order component, one obtains:
\begin{equation}
	\bfE \approx \sum_{k=1}^{m_{o+1}} \bfE^{(o+1,k)} \bfz^{\alphavec(o+1,k)} = \sum_{k=1}^{m_{o+1}} \bfR^{(o+1,k)} \bfz^{\alphavec(o+1,k)}.
\end{equation}
In particular, for the single master mode case, this expression can be simplified by employing the polar representation of the normal variables, yielding for the norm of the invariance error
\begin{equation}\label{eq:errorapprox01}
	\norm{\bfE} \approx \norm{\sum_{k=1}^{m_{o+1}} \bfR^{(o+1,k)} e^{i(\alphavec(p,k)_1-\alphavec(p,k)_2)\theta}} \left( \frac{\rho}{2} \right)^{(o+1)}.
\end{equation}
This equation approximates the numerator of~\cref{eq:inv_eq_error}, and requires much fewer computations. It could then already be used to simplify the criterion, instead of \cref{eq:norm_approx_HBNF}. It still depends, however, on the phase $\theta$ of the normal coordinates, and does not yield a direct analytic expression for the validity limit. To get rid of the angular dependence on $\theta$, and obtain an easy and rapid formulation, a first possibility is to compute~\cref{eq:errorapprox01} for $\theta=0$. Alternatively, it is possible to consider the exponential functions in \cref{eq:errorapprox01} as a projection basis for the vector $\bfE$, with $\bfR^{(o+1,k)}$ its coordinates in this base. Thus, introducing the residual vector $\bfR^{(p)}$ at order $p$, which gathers the residuals for all of the considered monomials:  
\begin{equation}
	\bfR^{(p)} = \left[ (\bfR^{(p,k_1)})^T \quad (\bfR^{(p,k_2)})^T \quad \ldots \quad (\bfR^{(p,m_p-1)})^T \quad (\bfR^{(p,m_p)})^T \right]^T,
\end{equation}
one can use this idea to compute an approximation of the norm of the invariance equation error as simply the norm of $\bfR^{(o+1)}$, which reads:
\begin{equation}
	\norm{\bfE} \approx \norm{\bfR^{(o+1)}} \left( \frac{\rho}{2} \right)^{(o+1)}.
\end{equation}
This last equation then immediately yields~\cref{eq:inv_eq_error_ANM} after substitution into \cref{eq:inv_eq_error}. Note that the idea of using only the first neglected term of the series has already been used in the context of the Asymptotic-Numerical Method (ANM), see {\em e.g.}~\cite{Cochelin1994}, leading to simplified formulations for determining the convergence range. Also, using the exponential functions as a projection basis is one of the core concepts of the Harmonic Balance Method, which inspired the proposed approximation. In a context more related to the parametrisation method, the same approach has also been used in~\cite{Grolet:HBNF}, where a very similar computation as the one used here with~\cref{eq:inv_eq_error_ANM}, is also proposed. More specifically, in~\cite{Grolet:HBNF}, a direct link between the Harmonic Balance and the use of complex normal form in the parametrisation method, is proposed, hence yielding a different style of parametrisation, which has been named as Harmonic Balance Normal Form (HBNF) style.

\section{Computation of the maximum value of the forcing amplitude for primary resonance} \label{sec:max_ampl_comp}

In this appendix, the relation between the amplitude of the forcing and of the normal variables for a primary resonance scenario is compared, in order to estimate {\em a priori} the maximum forcing that can be employed for a non-autonomous system. Given a system in primary resonance subjected to the parametrisation method with CNF style and reduced to a single master mode, its reduced dynamics can be written generically as \cite{Stabile:morfesym}
\begin{equation}
	\dot{z}_1 = f_1 z_1 + f_2 z_1^2z_2 + f_3 z_3,
\end{equation}
up to order 3 on the autonomous variables and order 1 on the non-autonomous ones, with the equation for $z_2$ being its complex conjugate. We decompose the coefficients into real and imaginary parts as $f_i = f_i^R + i f_i^I$, and in particular, we note that $f_3$ is always linear with respect to the forcing amplitude $\kappa$, and we write it as $f_3 = (c_3^R + ic_3^I) \kappa$. Then, introducing the polar representation of the normal variables, given by \cref{eq:polar_autonomous,eq:polar_non_autonomous}, separating real and imaginary parts and introducing the phase difference $\phi$ defined by \cref{eq:resonance} in order to render the equations autonomous, we find the following dynamical system:
\begin{subequations}
	\begin{align} 
		\dot{\rho} &= f_1^R \rho + f_2^R \frac{\rho^3}{4} + 2 f_3^R \cos{\phi} + 2 f_3^I \sin{\phi}  \label{eq:FRC_syst_a} \\
		\rho \dot{\phi} &= (f_1^I - \Omega) \rho + f_2^I \frac{\rho^3}{4} + 2 f_3^I \cos{\phi} - 2 f_3^R \sin{\phi} \label{eq:FRC_syst_b}.
	\end{align} \label{eq:FRC_syst}
\end{subequations}

Periodic orbits of the original system correspond to fixed points of \cref{eq:FRC_syst}. Then, if we equal the right-hand sides to zero and assume for the moment that $f_3^R \neq 0$, we can isolate $\cos{\phi}$ from \cref{eq:FRC_syst_a}:
\begin{equation} \label{eq:cos_int}
	\cos{\phi} = - \frac{f_1^R \rho + f_2^R \frac{\rho^3}{4} + 2 f_3^I \sin{\phi}}{2 f_3^R}.
\end{equation}
Substituting this expression into \cref{eq:FRC_syst_b} and defining $\alpha = \nicefrac{f_3^I}{f_3^R}$, we obtain, after algebraic manipulations
\begin{equation} \label{eq:sin_phi}
	\sin{\phi} = \frac{(f_1^I-\alpha f_1^R-\Omega)\rho + (f_2^I - \alpha f_2^R) \frac{\rho^3}{4}}{2 f_3^R (1+\alpha^2)},
\end{equation}
which can be substituted back into \cref{eq:cos_int}, yielding
\begin{equation} \label{eq:cos_phi}
	\cos{\phi} = -\frac{(f_1^R+\alpha f_1^I- \alpha \Omega)\rho + (f_2^R + \alpha f_2^I) \frac{\rho^3}{4}}{2 f_3^R (1+\alpha^2)}.
\end{equation}
These two expressions can be substituted into the identity $\sin^2{\phi} + \cos^2{\phi} = 1$, we get
\begin{equation}
	\left[(f_1^I-\alpha f_1^R-\Omega)\rho + (f_2^I - \alpha f_2^R) \frac{\rho^3}{4} \right]^2 + \left[ (f_1^R+\alpha f_1^I- \alpha \Omega)\rho + (f_2^R + \alpha f_2^I) \frac{\rho^3}{4} \right]^2 = 4 (f_3^R)^2 (1+\alpha^2)^2.
\end{equation}
By expanding the expressions on the left-hand side and passing the right-hand side to the left, it is possible to get a quadratic equation on $\Omega$ of the form
\begin{equation}
	a \Omega^2 + b \Omega + c = 0,
\end{equation}
with coefficients defined by
\begin{subequations}
	\begin{align}
		a &= \rho^2, \\
		b &= - 2f_1^I\rho^2 - f_2^I \frac{\rho^4}{2}, \\
		c &= |f_1|^2 \rho^2 + (f_1^I f_2^I + f_1^R f_2^R) \frac{\rho^4}{2} + |f_2|^2 \frac{\rho^6}{16} - 4 (f_3^R)^2 (1 + \alpha^2).
	\end{align}
\end{subequations}
The peak of the FRC is defined when this equation has a single solution, and thus we require its discriminant to be null, {\em i.e.} $b^2 - 4ac = 0$. This implies, by taking $\kappa$ out of $f_3$ and manipulating the expression:
\begin{equation}
	16 (c_3^R)^2 (1+\alpha^2) \kappa^2 = 4(f_1^R)^2 \rho^2 + 2f_1^R f_2^R \rho^4 + (f_2^R)^2 \frac{\rho^6}{4}.
\end{equation}
Therefore, we find the final expression for $\kappa$ as
\begin{equation}
	\kappa = \frac{2 f_1^R \rho + f_2 ^R \frac{\rho^3}{2}}{4c_3^R \sqrt{1 + \alpha^2}}.
\end{equation}
Although the above expression was deduced for the case $f_3^R \neq 0$, an analogous one can be computed by considering instead $f_3^I \neq 0$ (and one of these two conditions will always be verified). In this situation, introducing $\beta = \alpha^{-1} = \nicefrac{f_3^R}{f_3^I}$, the final expression for the forcing amplitude is
\begin{equation} \label{eq:max_kappa_fI}
	\kappa = \frac{2 f_1^R \rho + f_2 ^R \frac{\rho^3}{2}}{4c_3^I \sqrt{1 + \beta^2}}.
\end{equation}
In particular, this expression specializes to the one given for the Duffing oscillator in \cite{Stabile:morfesym}.

\section{Comparison of different strategies for the singularity criterion} \label{sec:sing_comp}

\begin{figure}[h]
	\centering
	\begin{subfigure}{0.49\textwidth}
		\includegraphics[width=\textwidth]{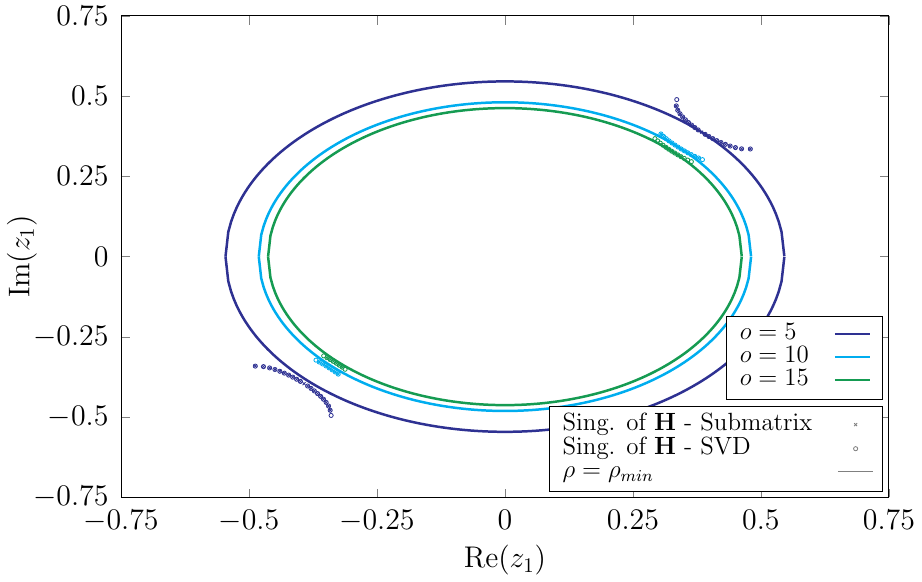}
		\caption{Undamped and unforced}
		\label{fig:Cantilever unforced - Singularity crtiterion - Comparison}
	\end{subfigure}
	\begin{subfigure}{0.49\textwidth}
		\includegraphics[width=\textwidth]{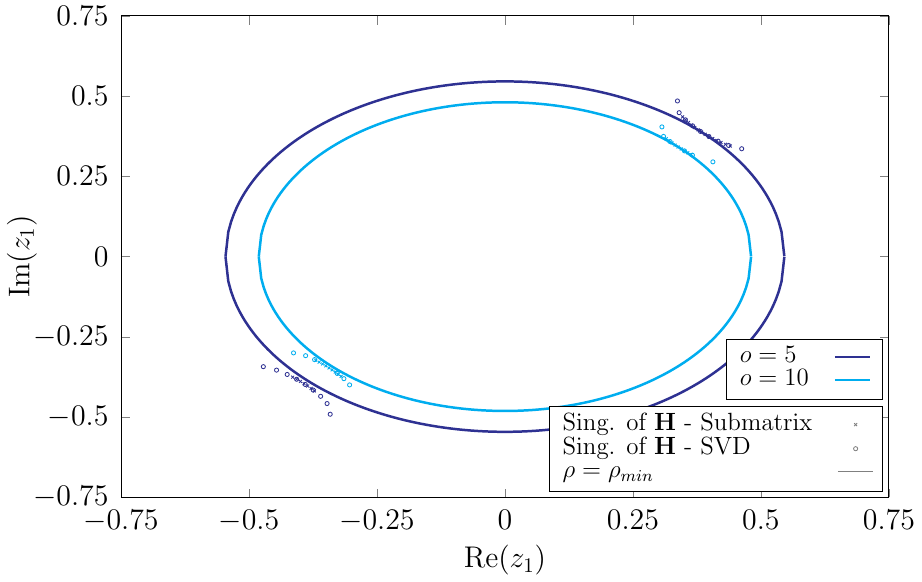}
		\caption{Primary resonance - $\kappa=7 \text{ m/s}^2$}
		\label{fig:Cantilever damped forced - Singularity crtiterion - Comparison}
	\end{subfigure}
	\caption{Validity limits calculated by the criterion of the singularity of the homological operator for the cantilever beam as a function of the parametrisation order. We compare the submatrix and SVD approaches. The lines represent circles of radius equal to the minimum amplitude among points of the same order.}
	\label{fig:Cantilever - Singularity crtiterion - Comparison}
\end{figure}

When employing the singularity criterion for non-square homological operators, as described in \cref{sec:Singularity}, a strategy needs to be chosen in order to determine the points where the operator's rank gets smaller. For the examples presented in the main text, we have chosen to calculate the determinant of a relevant submatrix and determine its zeroes by means of a root-finding algorithm. Another approach consists of considering the complete operator's non-null singular values, and, by means of a minimization approach, find the points where the smallest of them becomes zero. In this appendix, we compare these two approaches for the cantilever beam presented in \cref{sec:Beam}, both without damping and forcing and in primary resonance. The points where the operator's rank becomes smaller following the two approaches are presented in \cref{fig:Cantilever - Singularity crtiterion - Comparison}. We remark that the points computed following the two approaches are not exactly the same, but are located at the same zones and correspond to the same values of convergence radius. This, thus, validates the submatrix approach for this example. A situation where the SVD strategy could be preferred would be when no representative degrees of freedom (such as the displacement at the tip of the beam for this case) are readily available.

\end{document}